\documentclass[leqno,11pt,a4paper]{amsart}

\usepackage{jc}
\usepackage{tikz}
\usepackage{tikz-cd}
\usepackage[normalem]{ulem}
\usepackage[mathcal]{eucal}
\usepackage{comment}

\begin{document}

\author{Julian Chaidez}
\address{Department of Mathematics\\University of Southern California\\Los Angeles, CA\\90007\\USA}
\email{julian.chaidez@usc.edu}

\author{Yijie Pan}
\address{Department of Mathematics\\University of Southern California\\Los Angeles, CA\\90007\\USA}
\email{yijiepan@usc.edu}

\title[Pseudo-Anosov Reeb Flows And Contact Structures]{Pseudo-Anosov Reeb Flows And Contact Structures}

\begin{abstract} We introduce the notion of a pseudo-Anosov contact structure, which admits a type of singular contact form with pseudo-Anosov Reeb flow. We prove that contact homology detects the free homotopy classes of closed orbits of any pseudo-Anosov Reeb flow and that any pseudo-Anosov contact structure is universally tight and torsion free. Many applications are given, including new cases of the Finiteness Conjecture for transitive pseudo-Anosov flows. Our proofs use a flavor of contact homology graded by a free homotopy class of loops, defined for any contact manifold. We establish several properties of this type of contact homology that may be of independent interest.
\end{abstract}

\vspace*{-20pt}

\maketitle

\vspace*{-25pt}

\section{Introduction} \label{sec:introduction} 

A \emph{pseudo-Anosov flow} is a special type of flow on a closed 3-manifold that is Anosov away from a link of \emph{singular orbits}, each of which has a neighborhood modeled on the suspension of a standard map with a multi-pronged singular fixed point. Pseudo-Anosov flows are a generalization of Anosov flows that were originally introduced by Thurston (see e.g. \cite{thurston1997three}). They have since become a central topic at the intersection of low-dimensional topology and dynamics (cf. \cite{potrie2025anosov,barbot2013pseudo,barthelmé2022orbitequivalencespseudoanosovflows}).

\vspace{3pt}

Recently, a fruitful interface has emerged between the theory of pseudo-Anosov flows and symplectic topology. This is part of an emerging current of interactions between symplectic topology and hyperbolic dynamics more generally, leading to new developments in non-Weinstein Liouville geometry \cite{hozoori2024symplectic,hozoori2024regularity,hozoori2023anosovity,massoni2025anosov,cieliebak2022floer} and convex hypersurface theory \cite{chaidez2024robustly}. Of particular note, this interface has results in significant progress on the following fundamental conjecture.

\begin{conjecture*}[Finiteness] \label{conj:finiteness} There are only finitely many transitive pseudo-Anosov flows up to orbit equivalence on any closed 3-manifold. 
\end{conjecture*}

\noindent Recall that two flows are \emph{orbit equivalent} if there is a homeomorphism that maps the orbits of one flow to the orbits of the other. Moreover, a flow is \emph{transitive} if it possesses a dense orbit. Recently, Barthelmé-Bowden-Mann proved the finiteness of contact (or equivalently, Reeb) Anosov flows up to orbit equivalence, by combining tools from contact topology with major breakthroughs in the classification of pseudo-Anosov flows via the homotopy classes of their orbits, due to Barthelmé-Mann \cite{barthelme2024orbit} and Barthelmé-Frenkel-Mann \cite{barthelmé2022orbitequivalencespseudoanosovflows}. More recently, Zung obtained partial extensions of these results to other transitive pseudo-Anosov flows \cite{zung2024pseudoanosovrepresentativesstablehamiltonian} and Baldwin-Sivek-Zung \cite{baldwin2025pseudo} applied these methods to construct large families of distinct transitive pseudo-Anosov flows on hyperbolic L-spaces, addressing questions of Calegari, Min and Nonnino.

\vspace{3pt}

Building on these developments, this paper introduces the notion of a pseudo-Anosov contact structure on a 3-manifold, which admits a type of singular contact form whose Reeb flow is pseudo-Anosov. Many examples can be constructed by taking branched covers of contact manifolds with Anosov Reeb flows (e.g. those of Foulon-Hasselblatt \cite{foulon2013contact}). We show that contact homology can be used to detect the free homotopy classes of closed orbits for pseudo-Anosov Reeb flows (Theorem \ref{thm:contact_homology}). We also prove that pseudo-Anosov contact forms are universally tight (Corollary \ref{cor:tightness}) and torsion free (Theorem \ref{thm:torsion_free}). Our proofs use a flavor of contact homology graded by a free homotopy class of loops, defined for any contact manifold. We establish several properties of this type of contact homology that may be of independent interest.

\vspace{3pt}

As applications, we improve upon various results of Barthelm\'{e}-Bowden-Mann \cite{barthelme2024orbit}, Zung \cite{zung2024pseudoanosovrepresentativesstablehamiltonian} and Alves \cite{alves2016cylindrical}. We show that the pseudo-Anosov Reeb flow for a pseudo-Anosov contact structure is unique up to isotopy equivalence (Theorem \ref{thm:canonical_flow}), and that the pseudo-Anosov Reeb flow has the minimal set of orbit classes among all Reeb flows for the same contact structure (Theorem \ref{thm:orbit_classes}), addressing various questions posed recently by Barthelm\'{e} \cite{barthelme2025smorgaasbord}. We show that all Reeb flows for pseudo-Anosov contact structures have positive entropy (Theorem \ref{thm:entropy}). Finally, we prove Conjecture \ref{conj:finiteness} holds for the class of pseudo-Anosov Reeb flows (Theorem \ref{thm:finiteness}).

\subsection{Pseudo-Anosov Contact Forms And Smoothings.} \label{subsec:main_results} In order to state our main results more carefully, let us informally introduce the basic objects of study in this paper. 

\begin{definition*} \label{def:PA_Reeb_flow} A pseudo-Anosov flow $\Phi$ on a closed 3-manifold $Y$ is \emph{contact} (or alternatively \emph{Reeb}) if it is generated by the Reeb vector field $R$ of a contact form $\alpha$ with singularities in a link.
\end{definition*}

\noindent Here a \emph{contact form with singularities along a link} (Definition \ref{def:singularcontact}) is essentially a Lipschitz 1-form on a 3-manifold that is smooth and contact away from a singular link, and whose Reeb flow extends continuously across the link. Non-trivial examples of pseudo-Anosov Reeb flows can be produced by taking branched covers of pseudo-Anosov flows over links of closed orbits, as originally noted by Etnyre-Ghrist \cite{etnyre1998tightcontactstructures3manifolds}. In Section \ref{sec:pseudo-Anosov_contact_forms}, we will provide a brief review of pseudo-Anosov flows, followed by a more detailed discussion of the above definitions and constructions.

\vspace{3pt}

In Section \ref{sec:smoothings}, we describe a smoothing procedure for singular contact forms. This depends on a choice of smoothing data, but produces a contact structure that is independent of it. We call this the \emph{associated contact structure} to the singular contact form (or equivalently its Reeb flow).

\begin{definition*} \label{def:PA_contact_structure} A contact structure $\xi$ on a closed $3$-manifold $Y$ is \emph{pseudo-Anosov} if there is a pseudo-Anosov Reeb flow $\Phi$ with associated contact structure $\xi$.
\end{definition*}

\noindent By an appropriate choice of smoothing data, one may obtain a smoothed contact form whose Reeb flow closely mirrors the original pseudo-Anosov flow. In particular, since the closed orbits of a pseudo-Anosov flow are always non-contractible, which yields the following theorem.

\begin{theorem*}[Hypertightness] \label{thm:hypertight} A pseudo-Anosov contact structure $\xi$ is asymptotically hypertight.
\end{theorem*}

\noindent Any asymptotically hypertight contact manifold has non-vanishing full contact homology, as do its finite covers. Therefore, Theorem \ref{thm:hypertight} has the following corollary. 

\begin{corollary*}[Universal Tightness] \label{cor:tightness} Any pseudo-Anosov contact structure $\xi$ is universally tight.
\end{corollary*}

\subsection{Contact Homology In A Homotopy Class} \label{subsec:intro_CH} Our main result is a computation of a flavor of contact homology graded by free homotopy classes of loop, which we now introduce briefly. Given a contact manifold $(Y,\xi)$ with non-degenerate contact form, the contact homology 
\[
CH(Y,\xi;\Gamma) \qquad\text{in a homotopy class}\qquad \Gamma \in [S^1,Y]
\]
is the homology of a certain complex generated by both contractible Reeb orbits and Reeb orbits in the homotopy class $\Gamma$. The differential of a contractible orbit counts genus zero, SFT-type curves in the symplectization with one positive end at a contractible orbit and many negative orbits on contractible orbits. The differential of an orbit in the class $\Gamma$ counts similar curves, where the positive end and one of the negative ends must be at an orbit in the homotopy class $\Gamma$.

\begin{remark*} To our knowledge, while this flavor of contact homology is known to some experts (cf. Pardon \cite[\S 1.8]{pardon2019contact}), this is the first appearance of a detailed account. In the hypertight case, it specializes to cylindrical contact homology, which has been studied extensively \cite{bourgeois2002morse,bourgeois2005homologie,hutchings2017cylindrical,hutchings2022s,hutchings2014cylindrical,alves2016cylindrical}.\end{remark*}

 In Section \ref{sec:contact_homology}, we discuss this flavor of contact homology in detail. We also extend a number of useful structural results from cylindrical contact homology to this more general theory. Of particular note, we have the following bound for the contact homology using Giroux torsion. 

\begin{theorem*}[Torsion Bound] \label{thm:torsion_bound} Let $(Y,\xi)$ be an algebraically tight, closed contact 3-manifold with Giroux torsion $k$ in the homotopy class of an essential torus $T \subset Y$, and let $\Gamma \in \pi_1(T)$ be non-zero. Then either
\[
\on{rank}(CH(Y,\xi;\Gamma)) \ge 2k \qquad\text{or}\qquad \on{rank}(CH(Y,\xi;\Gamma^{-1})) \ge 2k
\] 
\end{theorem*}

\begin{remark*} Theorem \ref{thm:torsion_bound} generalizes a result of Bourgeois-Giroux \cite{bourgeois2005homologie}, who proved a similar inequality for cylindrical contact homology, assuming the existence of a hypertight contact form that restricts to the standard contact form on the embedded Giroux domain. The novelty of our result is the removal of these requirements, which is a major benefit of working with our chosen flavor of contact homology. Algebraic tightness is a minimal assumption since an overtwisted contact structure has vanishing contact homology and arbitrarily high Giroux torsion.
\end{remark*}

\begin{remark*}[Sutured Case] The flavor of contact homology constructed in this paper can also be constructed for convex sutured contact manifolds via the methods of Colin-Ghiggini-Honda-Hutchings \cite{cghh2011}. The torsion bound in Theorem \ref{thm:torsion_bound} extends to the sutured case by the same proof.
\end{remark*}

\subsection{Contact Homology And Giroux Torsion Of Reeb Pseudo-Anosovs} In Section \ref{sec:contact_homology_of_PA_structures}, we use the smoothing methods in Section \ref{sec:smoothings} to give the following (partial) computation of contact homology for pseudo-Anosov contact manifolds. A more precise statement is given in Theorem \ref{thm:body_contact_homology}.

\begin{theorem*}[Contact Homology] \label{thm:contact_homology} Let $\xi$ be a pseudo-Anosov contact structure. Then for any primitive free homology class $\Gamma$ and any pseudo-Anosov Reeb flow $\Phi$ inducing $\xi$, we have
\[
CH(Y,\xi;\Gamma) \neq 0 \qquad\text{if and only if}\qquad \Gamma \in \mathcal{P}(\Phi)
\]
\end{theorem*}

\begin{remark*} A number of calculations of contact homology and symplectic field theory in the pseudo-Anosov setting have appeared previously in works of Cotton-Clay \cite{cotton2009symplectic}, Colin-Honda \cite{colin2008reebvectorfieldsopen} and Zung \cite{zung2024pseudoanosovrepresentativesstablehamiltonian}. Our work has a number of advantages and novelties relative to these works. 

\vspace{3pt}

First, we utilize a flavor of contact homology whose foundations are fully established by the virtual fundamental cycle constructions of Pardon \cite{pardon2019contact}. This is in contrast to linearized contact homology as applied by Colin-Honda in \cite{colin2008reebvectorfieldsopen}, where an unestablished type of homotopy invariance is used \cite[Thm 3.2(2)]{colin2008reebvectorfieldsopen}, and rational symplectic field theory as used by Zung in \cite{zung2024pseudoanosovrepresentativesstablehamiltonian}. Second, we circumvent the need for the construction of specific (and delicate) smooth resolutions of pseudo-Anosov flows as used by Zung in \cite{zung2024pseudoanosovrepresentativesstablehamiltonian}. Instead, we develop a softer approach, whereby we construct a large class of smoothed contact forms suitable for computing contact homology, and then appeal to a stability property for the Lefschetz index near a singular orbit (see Section \ref{subsec:lefschetz_and_tracking}) in our computation of the contact homology chain groups. Finally, our result can be enhanced to incorporate the action filtration in a manner that is key for the applications to entropy (Theorem \ref{thm:entropy}). An analogous argument is seems unavailable in the rational SFT framework of \cite{zung2024pseudoanosovrepresentativesstablehamiltonian}.\end{remark*}

\noindent As a key application of Theorem \ref{thm:torsion_bound} and Theorem \ref{thm:contact_homology}, we can show that pseudo-Anosov contact structures are torsion free. See Section \ref{subsec:torsion_and_CH} for a review of terminology regarding Giroux torsion.

\begin{theorem*}[Torsion Free] \label{thm:torsion_free} A pseudo-Anosov contact structure $\xi$ has zero Giroux torsion.
\end{theorem*}

\begin{proof} Suppose otherwise and let $U = U_1 \subset Y$ be an embedded Giroux torsion domain with torsion one in the class of an essential torus (see Definition \ref{def:torsion}). Consider the subgroup
\[
\pi_1(U) \subset \pi_1(Y) \qquad\text{isomorphic to $\mathbb{Z}^2$}
\]
Fix a pseudo-Anosov Reeb flow $\Phi$ on $\xi$ and recall that closed orbits of $\Phi$ in a given homotopy class $\Gamma$ are in bijection with points in the associated bifoliated plane $P$ that are fixed by $\Gamma$ under the associated Anosov-like $\pi_1(Y)$-action on $P$ \cite[Ch 2]{barthelmé2025pseudoanosovflowsplaneapproach}. Moreover, by \cite[Cor 3.2.4]{barthelmé2025pseudoanosovflowsplaneapproach}, there are generators $\Gamma$ and $\Xi$ of $\pi_1(U)$ such that $\Xi$ acts freely on $P$. In particular, the primitive homotopy classes $\Xi$ and $\Xi^{-1}$ must contain no closed orbits. On the other hand, if $\xi$ has Giroux torsion then Theorem \ref{thm:giroux_torsion_CH} states that for any homotopy class $\Xi \in \pi_1(U)$, we have either
\[
CH(Y,\xi;\Xi) \neq 0 \qquad\text{or}\qquad CH(Y,\xi;\Xi^{-1}) \neq 0
\]
By Theorem \ref{thm:body_contact_homology}, this implies that the pseudo-Anosov Reeb flow $\Phi$ must have an orbit in one of these two homotopy classes. This is a contradiction. \end{proof}

\begin{remark*}[Torsion For Reeb Anosovs] Barthelm\'{e}-Bowden-Mann \cite{barthelme2024orbit} proved Theorem \ref{thm:torsion_free} in the Anosov setting (where no singular orbits are present) with a clever argument using Liouville-Anosov domains. Namely, any Anosov contact structure $\xi$ on a 3-manifold $Y$ is semi-fillable by a Liouville domain constructed from the Anosov flow \cite{hozoori2023anosovity,hozoori2024regularity,massoni2025anosov}, while Giroux torsion obstructs semi-fillability. In the pseudo-Anosov setting, no such argument is available. The proof in the general pseudo-Anosov case above was proposed to us by Zung in private communication.
\end{remark*}

As discussed in the introduction, Theorems \ref{thm:contact_homology} and \ref{thm:torsion_free} have a number of further applications. In the rest of this introduction, we state these results and provide most of their proofs.

\subsection{Canonical Reeb Flow} A fundamental result of Thurston \cite{thurston1988geometry} states that any two pseudo-Anosov representatives of a given mapping class of a surface are conjugate. In other words, any pseudo-Anosov map provides a canonical dynamical representative of its isotopy class. In \cite{zung2024pseudoanosovrepresentativesstablehamiltonian}, Zung established an analogue of this statement for transitive pseudo-Anosov flows on hyperbolic 3-manifolds, in terms of exact deformation classes of stable Hamiltonian structures. We prove the following analogue in the contact setting, with no hyperbolicity assumption.

\begin{theorem*}[Canonical Flow] \label{thm:canonical_flow} Let $\Phi$ and $\Psi$ be two pseudo-Anosov Reeb flows with the isotopic associated contact structures. Then $\Phi$ and $\Psi$ are isotopy equivalent.
\end{theorem*}

\noindent The proof of Theorem \ref{thm:canonical_flow} is analogous to the proof for Reeb Anosov flows in \cite{barthelme2024orbit}, and we present it here. It requires Theorem \ref{thm:contact_homology}, Theorem \ref{thm:canonical_flow}, several basic lemmas in Section \ref{sec:pseudo-Anosov_contact_forms} and the following remarkable classification result of Barthelem\'{e}-Frenkel-Mann \cite[Thm 1.1 and Prop 1.2]{barthelmé2022orbitequivalencespseudoanosovflows}.

\begin{theorem*}[Barthelm\'{e}-Frenkel-Mann]\label{thm:orbiteq} Let $\Phi,\Psi$ be two transitive pseudo-Anosov flows without a tree of scalloped regions. Then $\Phi$ and $\Psi$ are isotopy equivalent if and only if
\[\mathcal{P}(\Phi)=\mathcal{P}(\Psi)\]
Moreover, if there is no essential torus transverse to $\Phi$, then $\Phi$ does not have a tree of scalloped regions. 
\end{theorem*}

\begin{proof} (Theorem \ref{thm:canonical_flow}) Let $\Phi$ and $\Psi$ be two pseudo-Anosov Reeb flows with isotopic associated contact structures $\xi$ and $\eta$, respectively. By the isotopy invariance of contact homology
\[
CH(Y,\xi;\Gamma) \simeq CH(Y,\eta;\Gamma) \qquad\text{for any simple free homotopy class }\Gamma
\]
Theorem \ref{thm:contact_homology} then implies that a simple homotopy class $\Gamma$ contains an orbit of $\Phi$ if and only if it contains an orbit of $\Psi$. The set of orbit classes of a pseudo-Anosov flow is actually determined entirely by the simple orbit classes that is contains (c.f. Lemma \ref{lem:simple_vs_iterated_orbit_classes}). Therefore we know that $\mathcal{P}(\Phi) = \mathcal{P}(\Psi)$. Finally, the pseudo-Anosov Reeb flows $\Phi$ and $\Psi$ are transitive (Lemma \ref{lem:transitive}) and admit no transverse embedded tori (Lemma \ref{lem:no_transversals}). The result thus follows from Theorem \ref{thm:orbiteq}. 
\end{proof}

\subsection{Orbit Classes} In the recent survey \cite{barthelme2025smorgaasbord}, Barthelm\'{e} posed several questions about canonical Reeb dynamics. Generally, Barthelm\'{e} asked how related (dynamically) two Reeb flows for the same contact structure must be \cite[Question 1]{barthelme2025smorgaasbord}. More precisely, he posed the following question.

\begin{question*} \cite{barthelme2025smorgaasbord} \label{qu:Barthelme_minimal_flow} Let $\xi$ be a contact structure. Under what conditions on $\xi$ does there exist a Reeb flow $\Phi_{\on{min}}$ for $\xi$ such that, for any other Reeb flow $\Phi$ for $\xi$, we have $\mathcal{P}(\Phi_{\on{min}}) \subset \mathcal{P}(\Phi)$?
\end{question*}

\noindent A basic property of our variant of contact homology is that it vanishes in a particular homotopy class if there is a contact form with no orbits in that class (Lemma \ref{lem:orbit_existence}). Thus, as a corollary of Theorem \ref{thm:contact_homology}, we immediately acquire the following answer to Question \ref{qu:Barthelme_minimal_flow}.

\begin{theorem*}[Orbit Classes] \label{thm:orbit_classes} Let $\Phi_{PA}$ be a pseudo-Anosov Reeb flow for a contact structure $\xi$. Then for any smooth Reeb flow $\Phi$ for $\xi$, there is an inclusion
\[
\mathcal{P}(\Phi_{PA}) \subseteq \mathcal{P}(\Phi)
\]
\end{theorem*}

\noindent 

\subsection{Entropy} Several recent works \cite{alves2016cylindrical,cineli2025barcode,alves2021braid,alves2019dynamically,abbondandolo2023entropy,fernandes2024barcode,alves2023c,fender2023barcode,alves2022reeb} have focused on the close relationship between topological entropy and Floer homology for symplectomorphisms and Reeb flows. 

\vspace{3pt}

In \cite{alves2016cylindrical}, Alves proved the first result in this direction, showing that the existence of a hypertight Reeb flow with exponential homotopical growth in cylindrical contact homology forces all Reeb flows for the same contact structure to have positive entropy. Anosov Reeb flow are the basic example of such hypertight Reeb flows. We generalize this result to pseudo-Anosov Reeb flows.

\begin{theorem*}[Entropy] \label{thm:entropy} Any Reeb flow $\Phi$ for a pseudo-Anosov contact structure $\xi$ has positive entropy. 
\end{theorem*}

\noindent This result follows from a lower bound for the entropy in terms of the homotopical growth rate of the flavor of contact homology discussed in Section \ref{sec:contact_homology}. We give a formal proof in Section \ref{sec:contact_homology_of_PA_structures}.

\subsection{Finiteness} Finally, our work yields a new special case of Conjecture \ref{conj:finiteness}, namely the finiteness of pseudo-Anosov Reeb flows up to isotopy equivalence. This generalizes the finiteness results of Barthelm\'{e}-Bowden-Mann \cite{barthelme2024orbit} for Anosov Reeb flows and the finiteness result of Zung \cite{zung2024pseudoanosovrepresentativesstablehamiltonian} for transitive pseudo-Anosov flows on homology spheres admitting a positive Birkhoff section   . 

\begin{theorem*}[Finiteness] \label{thm:finiteness} There are only finitely many pseudo-Anosov Reeb flows up to isotopy equivalence on any closed orientable $3$-manifold $Y$.
\end{theorem*}

\begin{proof} It follows from Theorem \ref{thm:canonical_flow} and Corollary \ref{cor:tightness} that the set of pseudo-Anosov Reeb flows up to isotopy equivalence on $Y$ is smaller than (or the same size as) the set of torsion free contact structures up to isotopy on $Y$. On the other hand, a closed orientable 3-manifold has finitely many torsion free contact structures up to isotopy by Colin-Giroux-Honda \cite{cgh2003}.
\end{proof}

\subsection{Questions And Open Problems.} We conclude this introduction with a discussion of several open problems related to this work. Some will be the subject of future work by the authors.

\vspace{3pt}

First, note that some closed 3-manifolds do not admit a pseudo-Anosov contact structure. Indeed, there are 3-manifolds that do not even admit a tight contact structure (cf. Honda-Etnyre \cite{e2001}) and many 3-manifolds (including hyperbolic 3-manifolds) with no pseudo-Anosov flows (e.g. the Weeks manifold). Thus the following question seems quite natural and interesting.

\begin{question*}[Existence] Which closed 3-manifolds $Y$ admit a pseudo-Anosov contact structure?
\end{question*}

Next, recall that Marty \cite{marty2023skewed} recently proved that all skew $\R$-covered Anosov flows in dimension three are orbit equivalent to a Reeb flow. This confirmed a longstanding conjecture on the equivalence of the Reeb and skew $\R$-covered conditions (c.f. Barbot \cite{barbot2001plane}) and gives a completely dynamical criterion for Reebness via the bifoliated plane. There is also an alternative characterization of skew $\R$-coveredness due to Asaoka-Bonatti-Marty \cite{asaoka2024oriented} via the existence of an oriented Birkhoff section with positive boundary.

\begin{question*}[Reeb Criteria] Let $\Phi$ be a pseudo-Anosov flow. Is there a necessary and sufficient criterion via the bifoliated plane (or via Birkhoff sections) that determines if $\Phi$ is Reeb?
\end{question*}

Finally, previous works of Colin-Honda \cite{colin2008reebvectorfieldsopen} studied contact manifolds admitting open books with pseudo-Anosov monodromy. In particular, \cite{colin2008reebvectorfieldsopen} established the exponential action growth rate of a certain linearized contact homology group under the hypothesis that the fractional Dehn twist coefficient $c = k/n$ satisfies $k \ge 3$. It is interesting to ask about the relationship between the methods of \cite{colin2008reebvectorfieldsopen} and this paper. Precisely, we pose the following question.

\begin{question*} Is there a characterization of pseudo-Anosov contact manifolds in terms of the fractional Dehn twist coefficient of a supporting pseudo-Anosov open book?
\end{question*}

\subsection*{Acknowledgements} We would like to thank Thomas Barthelm\'{e}, Audrey Rosevear, and especially Jonathan Zung for proposing the approach to Theorem \ref{thm:torsion_free}. JC was partially supported by NSF DMS-2446019 and by the US-Israel Binational Science Foundation under award 2024157.

\newpage

\section{Pseudo-Anosov Contact Manifolds} \label{sec:pseudo-Anosov_contact_forms} 

In this section, we formally introduce pseudo-Anosov contact forms and Reeb flows. 

\subsection{Pseudo-Anosov Flows} \label{subsec:pseudo-Anosov_flows} We start with a brief review of the theory of pseudo-Anosov maps and flows. Useful references for this material include Barthelm\'{e}-Mann \cite{barthelmé2025pseudoanosovflowsplaneapproach} and Tsang \cite{tsang2024horizontalgoodmansurgeryequivalence}.

\vspace{3pt}

We start by recalling the standard model for a singular orbit of a pseudo-Anosov. Consider the map of the plane for any $\lambda > 0$ given by
\[
A_\lambda:\R^2 \to \R^2 \qquad\text{given by}\qquad A_\lambda(x,y) = (\lambda x,\lambda^{-1}y)
\]
Also consider the standard branched covering map of the plane for any integer $n \ge 2$ given by
\[
\pi_n:\R^2 \to \R^2 \qquad\text{given in radial coordinates by}\qquad \pi_n(r,\theta) = (r,n\theta/2)
\]
\begin{definition}[Standard PA Map] The \emph{standard pseudo-Anosov map} with stretching factor $\lambda$, $n$ prongs and rotation $k$ is the unique map
\[
\phi_{n,k,\lambda}:\R^2 \to \R^2
\]
that lifts the map $A_\lambda$ via the branched cover $\pi_n$ and sends the ray $x$-axis $L = 0 \times \R_+$ in the plane to the rotation of $L$ by $2\pi k/n$ radians. The standard \emph{stable} and \emph{unstable} foliations 
\[
F_n^s \qquad\text{and}\qquad F^u_s
\]
are the foliations of $\R^2$ away from the origin that respectively lift the foliations of $\R^2$ by vertical and horizontal lines. Note that the standard map $\phi_{n,k,\lambda}$ preserves the foliations $F^s$ and $F^u$.\end{definition}

\begin{definition}[Standard PA Flow] The \emph{standard pseudo-Anosov flow}  with stretching factor $\lambda$, $n$ prongs and rotation $k$ is the suspension flow
\[
\Phi_{n,k,\lambda} \qquad\text{on the mapping torus}\qquad U_{n,k,\lambda} 
\]
Here recall that the mapping torus of a map is the quotient of $\R \times \R^2$ by the equivalence relation that identifies $(t+1,x)$ with $(t,\phi_{n,k,\lambda}(x))$. The standard \emph{weak stable foliation} $W^s_{\on{std}}$ and \emph{weak unstable foliation} $W^u_{\on{std}}$ are the unique foliations that descend from the foliations
\[
\R \times F^s_n \qquad\text{and}\qquad \R \times F^u_n \qquad\text{on}\qquad \R \times \R^2
\]\end{definition}

\begin{definition}[Pseudo-Anosov Flow]\label{pseudo-anosov flow}
A \emph{pseudo-Anosov flow} $\Phi$ on a closed 3-manifold $Y$ is a continuous flow such that
\[\Phi \text{ is smooth and Anosov away from a smoothly embedded link }\Gamma \subset Y\]
and such that each component $\gamma$ of $\Gamma$ admits a smooth, bi-Lipschitz tubular neighborhood chart
\[
\iota:\R/\Z \times \D \to Y \qquad\text{restricting to a diffeomorphism}\qquad \R/\Z \times 0 \to \gamma
\]
that sends orbits of $\Phi_{n,k,\lambda}$ to orbits of $\Phi$ and that sends the weak stable and unstable foliation $W^s_{\on{std}}$ and $W^u_{\on{std}}$ to the corresponding weak stable and unstable foliations of the Anosov flow $\Phi$.
\end{definition}

\begin{remark}[Pseudo-Anosov Definition] Our notion of pseudo-Anosov flow is essentially equivalent to the one given in Tsang \cite[Def 2.3]{tsang2024horizontalgoodmansurgeryequivalence}, with an additional assumption that the chart $\iota$ is smooth along $\Gamma$. This latter condition can always be achieved by simply changing the smooth structure on $Y$ near $\Gamma$ by asserting that $\iota$ is a chart map. Note that this does not change the diffeomorphism type of $Y$. More generally, a flow satisfying Definition \ref{pseudo-anosov flow} is often called a \emph{smooth pseudo-Anosov flow}, in contrast to the purely topological definition appearing in e.g. \cite{barthelmé2025pseudoanosovflowsplaneapproach}. In \cite{Agol_2024}, the authors prove that they are equivalent in the sense that any topological pseudo-Anosov flow is orbit equivalent to a smooth one. Thus we will not make any such distinction in this paper.

\vspace{3pt}

Also note that our definition does not allow 1-pronged singularities. In the smooth (non-contact) setting, pseudo-Anosov flows with 1-pronged singularities exhibit much more flexibility than ordinary pseudo-Anosov flows, and we expect a similar phenomenon in the contact case. 
\end{remark}

We will require a few basic results about pseudo-Anosov flows. 

\begin{lemma} The periodic orbits of a pseudo-Anosov flows are isolated.
\end{lemma}

\begin{proof} This holds for the smooth orbits since they are all hyperbolic (and thus non-degenerate). It holds for the singular orbits since $0$ is the unique fixed point of the standard map for any set of parameters (and thus the corresponding closed orbit is isolated in the suspension).
\end{proof}

\noindent We will also need the following elementary facts about the homotopy classes of closed orbits of a pseudo-Anosov flow. 

\begin{lemma}\cite[Lem 4.2]{zung2024pseudoanosovrepresentativesstablehamiltonian} \label{lem:simple_vs_iterated_orbit_classes}
    Let $\Phi$ be a transitive pseudo-Anosov flow with a closed orbit in the homotopy class $\Gamma^k$ for a primitive class $\Gamma \in \pi_1(Y)$. Then there is a closed orbit in either the homotopy class $\Gamma$ or $\Gamma^{-1}$.
\end{lemma}

\begin{lemma}\cite[Lem 4.3]{zung2024pseudoanosovrepresentativesstablehamiltonian} \label{lem:orbit_types_in_homotopy_classes}
    Let $\Gamma$ be a primitive homotopy class containing a periodic orbit of a transitive pseudo-Anosov flow $\Phi$. Then $\Gamma$ is represented by
    \begin{enumerate}
        \item a unique negative hyperbolic orbit, or 

        \item a unique rotating hyperbolic orbit, or

        \item (possibly distinct) positive hyperbolic or non-rotating singular orbits. 
        
    \end{enumerate}
\end{lemma}

\begin{lemma}\cite[Prop 1.3.6]{barthelmé2025pseudoanosovflowsplaneapproach} \label{lem:non_contractible_orbit} Let $\eta$ be a closed orbit of a transitive pseudo-Anosov flow $\Phi$ on a closed 3-manifold $Y$. Then $\eta$ is non-contractible.
\end{lemma}

\subsection{Singular And Pseudo-Anosov Contact Forms} We next introduce classes of singular contact forms that will play a key role in this paper. Fix a closed oriented 3-manifold
\[
Y \qquad\text{with a smooth embedded link}\qquad \Gamma \subset Y
\]

\begin{definition}[Singular Contact Form]\label{def:singularcontact} A \emph{contact form $\alpha$ on $Y$ with singularities on $\Gamma$} is a Lipschitz 1-form $\alpha$ such that
\begin{itemize}
    \item the 1-form $\alpha$ is smooth and contact on the complement of $\Gamma$
    \vspace{3pt}
    \item the Reeb vector field on the complement of $\Gamma$ extends to a continuous vector field on $Y$ that is tangent to $\Gamma$.
\end{itemize}
The \emph{Reeb vector field} $R$ of a contact form $\alpha$ with singularities on a link $\Gamma$ is the unique extension of the Reeb vector field on the complement of $\Gamma$. 
\end{definition} 

\begin{definition}[Pseudo-Anosov Contact Form] A contact form $\alpha$ with singularities in $\Gamma$ is \emph{pseudo-Anosov} if the associated Reeb vector field $R$ generates a pseudo-Anosov flow $\Phi$. 

\vspace{3pt}

Conversely, a pseudo-Anosov flow $\Phi$ is called \emph{Reeb} (or equivalently \emph{contact}) if $\Phi$ is the Reeb flow of a contact form with singularities in a link.
\end{definition}

\begin{example}[Contact/Reeb Anosov Flows] Any Anosov Reeb flow $\Phi$ is a pseudo-Anosov Reeb flow, where the corresponding singular contact form has empty singular link. The most basic examples are the cogeodesic flows on a closed hyperbolic surface, and more general examples were constructed by Foulon-Hasselblatt \cite{foulon2013contact}. 
\end{example}

\begin{example}[Branched Covers] More generally, one can construct a large class of pseudo-Anosov (and not simply Anosov) Reeb flows by the following branched cover construction (cf. Etnyre-Ghrist \cite{etnyre1998tightcontactstructures3manifolds}). Fix a pseudo-Anosov (or simply Anosov) Reeb flow
\[\Phi:\R \times Y \to Y \qquad\text{for a pseudo-Anosov contact form $\alpha$ on $Y$}\]
Denote the link of singularities of the contact form by $\Gamma$. Fix an additional link of Reeb orbits
\[\Xi \subset Y\]
Take any branched cover $\pi:\tilde{Y} \to Y$ that branches cover the link $\Gamma$ and let $\widetilde{\Gamma}$ and $\widetilde{\Xi}$ denote the links in the cover that projects to $\Gamma$ and $\Xi$ respectively. Then the 1-form
\[
\pi^*\alpha \qquad\text{on the covering 3-manifold}\qquad \widetilde{Y}
\]
is a pseudo-Anosov contact form with singular link $\widetilde{\Gamma} \cup \widetilde{\Xi}$ in the sense of Definition \ref{def:singularcontact}. Indeed, the pullback $\pi^*\alpha$ is smooth away from $\widetilde{\Gamma}$ and contact away from $\widetilde{\Gamma} \cup \widetilde{\Xi}$. The Reeb vector field $\widetilde{R}$ of $\pi^*\alpha$ is the lift of the Reeb vector field of $\alpha$ away from $\widetilde{\Xi}$ and it extends continuously to $\widetilde{\Xi}$ such that $\pi_*\widetilde{R} = R$ pointwise. Finally, we note that the lift of a pseudo-Anosov flow by a branched cover is well-known to be pseudo-Anosov (cf. Fenley \cite{fenley2003pseudo}).  \end{example}

A key property of pseudo-Anosov Reeb flows is that they have no transverse closed surfaces (and thus no scalloped regions cf. \cite{barthelmé2022orbitequivalencespseudoanosovflows}). We next prove this important fact.

\begin{lemma}[No Transversals] \label{lem:no_transversals} Fix a pseudo-Anosov Reeb flow $\Phi$ on a closed oriented 3-manifold $Y$ and let $\Sigma \subset Y$ be an embedded closed surface. Then $\Phi$ is not transverse to $\Sigma$. 
\end{lemma}

\begin{proof} Let $R$ be the Reeb vector field of the singular contact form $\alpha$ that generates $\Phi$, and let $\Gamma \subset Y$ denote the link of singular orbits. Also fix an auxilliary metric on $Y$.

\vspace{3pt}

Assume that $R$ is transverse to a closed embedded surface $\Sigma \subset Y$. Then $\Gamma$ is transverse to $\Sigma$, and thus intersects $\Sigma$ at a set of finite points $P \subset \Sigma$. We can decompose $\Sigma$ as
\[\Sigma = \Sigma_\epsilon \cup U_\epsilon \qquad\text{for small}\qquad \epsilon > 0\]
where $U_\epsilon$ is the union of the open $\epsilon$-balls around $P$ with respect to the auxilliary metric and $S_\epsilon$ is the complement of $U_\epsilon$. Since $S_\epsilon$ is contained in the region of $Y$ where $\alpha$ is smooth and the Reeb vector field $R$ is transverse to $\Sigma$, we have
\begin{equation} \label{eq:no_transversals}
d\alpha|_{S_\epsilon} > 0 \qquad\text{and thus}\qquad \int_{S_\epsilon} d\alpha = \int_{\partial S_\epsilon} \alpha > 0
\end{equation}
Note that $S_\delta \subset S_\epsilon$ if $\delta > \epsilon$, and so the integral (\ref{eq:no_transversals}) must stay positive as $\epsilon \to 0$. On the other hand, the length of the boundary of $\partial S_\epsilon = \partial U_\epsilon$ with respect to the auxilliary metric goes to zero as $\epsilon \to 0$. Since $\alpha$ is continuous, the integral (\ref{eq:no_transversals}) goes to zero as $\epsilon \to 0$. This is a contradiction.\end{proof}

\begin{corollary}[Transitive]\label{lem:transitive} A pseudo-Anosov Reeb flow $\Phi$ on a closed oriented 3-manifold $Y$ is transitive.
\end{corollary}
\begin{proof}
    Assume $\Phi$ is non-transitive, then there is a transverse torus $\Sigma$ to $\Phi$ by \cite[Proposition 2.7]{Mosher92} or \cite[Proposition 2.38]{barthelmé2022orbitequivalencespseudoanosovflows}, which contradicts Lemma \ref{lem:no_transversals}. 
\end{proof}

\section{Smoothings} \label{sec:smoothings}

In this section, we describe a smoothing construction for contact forms with singularities along a link and we analyze the properties of the smoothed Reeb flows. 

\subsection{Whitney Differential Forms} We start by recalling the notion of a Whitney differential form (cf. Whitney \cite{whitney2012geometric}), which is a differential form that admits a continuous exterior derivative. 

\begin{definition}[Whitney Forms] \label{def:whitney} A continuous $k$-form $\theta$ on a manifold $M$ is called \emph{Whitney} if the weak exterior derivative $d\theta$ is continuous. Equivalently, there is a continuous $(k+1)$-form
\[
\omega \qquad\text{such that}\qquad \int_\sigma \omega = \int_{\partial \sigma} \theta \text{ for all smooth simplices $\sigma$}
\]
\end{definition}

\begin{definition}[Whitney Topology] \label{def:whitney_topology} The \emph{Whitney topology} on the space of Whitney differential forms is induced by the $C^0$-topology on the differential form and the exterior derivative.
\[
\theta_i \xrightarrow{\text{Whitney}} \theta \qquad\text{if}\qquad (\theta_i,d\theta_i) \xrightarrow{C^0} (\theta,d\theta)
\]
\end{definition}

\begin{lemma} \label{lem:smooth_dense_Whitney} \cite[Lem 16(d)]{whitney2012geometric} Smooth differential forms are dense in the Whitney topology.
\end{lemma}

The notion of a contact form generalizes directly the setting of Whitney forms. We will need to use this generalization, along with a few elementary properties.

\begin{definition}[Whitney Contact Form] A \emph{Whitney contact form} $\alpha$ on a manifold $Y$ of dimension $2n+1$ is a Whitney 1-form such that the continuous differential form
\[\text{$\alpha \wedge d\alpha^n$ is nowhere vanishing.}\]
The \emph{Reeb vector field} $R$ is the unique continuous vector field such that $\iota_Rd\alpha = 0$ and $\iota_R\alpha = 1$. \end{definition}

\begin{lemma} \label{lem:whitney_contact_structure} Any Whitney contact form $\alpha$ on a compact manifold $Y$ determines a unique smooth contact structure $\xi$ up to isotopy (and thus up to contactomorphism).
\end{lemma}

\begin{proof} The set of smooth differential forms $\theta$ with $\|\theta - \alpha\|_{C^0} \le \epsilon$ and $\|d\theta - d\alpha\|_{C^0} \le \epsilon$ is convex and non-empty by Lemma \ref{lem:smooth_dense_Whitney}. Moreover, for sufficiently small $\epsilon$, each element of this set is contact since $\|\theta \wedge d\theta^n - \alpha \wedge d\alpha^n\| \le C \epsilon$ for appropriate constant $C$. This proves the result.  \end{proof}

\subsection{Construction Of Smoothings}
We next construct the appropriate class of Whitney (and smooth) contact forms that approximate a singular contact form. Fix a closed 3-manifold
\[
Y \qquad\text{with a contact form $\alpha$ with singularities along a link}\qquad \Gamma \subset Y
\] 

\begin{definition}[Smoothing Chart] \label{def:smoothing_chart} A \emph{smoothing chart} $(U,\iota,g)$ for the singular contact form $\alpha$ is a choice of a tubular neighborhood chart and a smooth function
\[
\iota:\Gamma \times \mathbb{D} \overset{\simeq}{\longrightarrow} U \subset Y \qquad\text{and}\qquad g:[0,1]_r \to \R
\]
such that $\iota$ restricts to the inclusion $0 \times \Gamma \to \Gamma$ along $0 \times \Gamma$, $g$ is a diffeomorphism on $(0,1)$, the derivatives of $g$ of all orders vanish at $0$ and $g(r)=r$ around $1$.
\end{definition}

\begin{remark} Given a choice of smoothing chart $(U,\iota,g)$, there is a smooth homeomorphism
\begin{equation} \label{eq:rho_map}
\rho:Y \to Y
\end{equation}
defined as the identity outside of $U$ and by $\rho(t,r,\theta) = (t,g(r),\theta)$ within $\Gamma \times \mathbb{D} \simeq U$. This map depends on the smoothing chart, but we will suppress that dependence notationally. \end{remark}

\begin{lemma}\label{lem:pullback_is_smooth} The pullback $\rho^*\alpha$ of $\alpha$ by the map $\rho:Y \to Y$ associated to a smoothing chart $(U,\iota,g)$ is a Whitney differential form on $Y$ such that
\[
\rho^*\alpha \text{ is smooth and contact on $Y \setminus \Gamma$}\qquad\text{and}\qquad d(\rho^*\alpha)\text{ vanishes along $\Gamma$}
\]  
\end{lemma}
\begin{proof} Since $\rho$ is smooth and $\alpha$ is smooth and contact away from $\Gamma$, the pullback $\rho^*\alpha$ is continuous, smooth away from $\Gamma$ and contact away from $\Gamma$. It suffices to show that the exterior derivative
\[
d(\rho^*\alpha) = \rho^*(d\alpha) \qquad\text{defined on}\qquad Y \setminus \Gamma
\]
extends continuously as $0$ across $\Gamma$. To prove this, assume without loss of generality that $\Gamma$ has a single component. By Definition \ref{def:singularcontact}, we may fix a parametrization $\R/\Z_t \simeq \Gamma$ such that $R|_\Gamma$ is a positive multiple of $\partial_t$. In the coordinates $\R/\Z \times \mathbb{D} \simeq U$, we may write
\[
d\alpha = dt \wedge \theta + \omega \qquad\text{and}\qquad R = \partial_t + V
\]
where $\theta$ is a smooth $\R/\Z$-dependent family of 1-forms on $\D \setminus 0$, $\omega$ is a smooth $\R/\Z$-dependent family of 1-forms on $\D \setminus 0$, and $V$ is a Lipschitz vector field on $\R/\Z \times \D$ that vanishes along $\R/\Z \times 0$. Moreover, since $\alpha$ is Lipschitz, $\theta$ and $\omega$ are bounded. Now note that
\[
|\theta| = |\iota_{\partial_t}(d\alpha)| \le |\iota_R(d\alpha)| + |\iota_Vd\alpha| \le |V| \cdot |d\alpha|
\]
Thus $|\theta|$ and $|\rho^*\theta|$ both limit to zero near $0 \times \Gamma$. Similarly, we have
\[
|\rho^*\omega| \le |\on{det}(d\rho)| \cdot |\omega|\le |g'(r)| \cdot |d\alpha|
\]
which limits to zero at $0 \times \Gamma$. This proves that $|\rho^*(d\alpha)|$ limits to zero at $0 \times \Gamma$. \end{proof}

\begin{definition}[Smoothing Function]\label{def:smoothing_function} A \emph{smoothing function} for a singular contact form $\alpha$ with respect to a smoothing chart $(U,\iota,g)$ is a smooth non-negative function
\[
\chi:[0,1]_r \to [0,\infty)
\]
such that there is an $\epsilon > 0$ and a constant $A > 0$ with
\begin{equation} \label{eq:smoothing_function_1}
 \text{$\chi(r) = Ar^2$ for $r < \epsilon$}\qquad\text{and}\qquad \text{$\chi(r) = 0$ for $r > 1-\epsilon$} 
\end{equation}
and such that the following Whitney 1-form given by the formula
\begin{equation}\label{eq:definition_of_alpha_chi}
    \alpha_\chi|_U = \rho^*\alpha+\chi(r)d\theta \qquad\text{and}\qquad \alpha_\chi|_{Y \setminus U} = \alpha
\end{equation}
is a contact form (see Definition \ref{def:whitney}). We denote the space of smoothing functions by
\[\mathcal{X} = \mathcal{X}(U,\iota,g)\]
\end{definition}

\begin{remark}[Coordinate Formulas] \label{rmk:coordinate_formulas} It will be useful to have formulas for the volume form and Reeb vector field for a smoothing function in the tubular coordinates. Precisely, let
\[
G:U \to \R \qquad\text{be the unique function with}\qquad \rho^*\alpha \wedge \rho^*d\alpha = G \cdot dt \wedge dr \wedge d\theta
\]
Given a fixed smoothing function $\chi$, the associated volume form is given by
\begin{equation} \label{eq:volume_form}
\alpha_\chi \wedge d\alpha_\chi = (G + H_\chi) \cdot dt \wedge dr \wedge d\theta \qquad\text{where}\qquad H_\chi = \chi'(r) \cdot \rho^*\alpha(\partial_t) + \chi(r) \cdot \rho^*d\alpha(\partial_t,\partial_r)
\end{equation} 
Similarly, the Reeb vector field $R_\chi$ associated toassociated to the smoothing function $\chi$ is given by
\begin{equation} \label{eq:reeb_vector_field}
R_\chi = \frac{1}{G + H_\chi} \cdot \big(F_t\partial_t + F_r\partial_r + F_\theta\partial_\theta\big)
\end{equation}
where the functions $F_\bullet$ are given by the formulas
\[
F_t = \rho^*d\alpha(\partial_r,\partial_\theta) + \chi'(r) \qquad F_r = \rho^*d\alpha(\partial_\theta,\partial_t) \quad\text{and}\quad F_\theta = \rho^*d\alpha(\partial_t,\partial_r)
\]
It is easily checked that $R_\chi$ satisfies the identity $\iota_{R_\chi}(\alpha_\chi \wedge d\alpha_\chi) = d\alpha_\chi$, which uniquely characterizes the Reeb vector field of a contact form. \end{remark}

\begin{lemma}[Existence] \label{lem:the_set_of_smoothing_functions} The space $\mathcal{X}$ of smoothing functions is non-empty and contractible. More precisely, fix a compact smooth manifold $S$ and a continuous family of smooth functions   
\[
\chi_s:[0,1]_r \to [0,\infty) \qquad\text{parametrized by $s \in S$ and satisfying (\ref{eq:smoothing_function_1})}
\]
Then for sufficiently small $\epsilon > 0$, the Whitney 1-forms $\alpha_{\epsilon \chi_s}$ are contact for each $s \in S$.
\end{lemma}
\begin{proof} Note that $\rho^*\alpha \wedge \rho^*d\alpha$ is a smooth volume form away from $\Gamma$. Since $\alpha_{\epsilon\chi} \to \alpha$ in $C^\infty$ as $\epsilon \to 0$ away from $\Gamma$, the lemma holds outside of $U$. Thus it suffices to show that $\alpha_{\epsilon\chi_s}$ is contact in the tubular coordinates chart $U$ for small $\epsilon$. By Remark \ref{rmk:coordinate_formulas} and the formula (\ref{eq:volume_form}), this is equivalent to showing that there are constants $C,\delta > 0$ and a small $\epsilon$ such that
\[
G + H_{\epsilon \chi_s} = G + \epsilon H_{\chi_s} \ge Cr \text{ for $r < \delta$} \qquad\text{and}\qquad G + \epsilon H_{\chi_s} > 0 \text{ for $r > 0$}
\]
for all parameters $s$. Note that the Reeb vector field $R$ of $\alpha$ is parallel to $\partial_t$ along $\Gamma$, the contact form $d\alpha$ is Lipschitz and $\rho^*\alpha \wedge d\alpha$ is a volume form away from $\Gamma$. Thus we know that
\[
\rho^*\alpha(\partial_t) > 0 \text{ near $r = 0$} \qquad\rho^*d\alpha(\partial_t,\partial_r)\text{ is bounded} \qquad\text{and}\qquad G \ge 0
\]
Moreover, by assumption $\chi_s(r) = A_sr^2$ near $r = 0$  for a parameter $A_s > 0$ depending continuously on $s$ and $G \ge 0$. Thus there is a small $\delta > 0$ independent of $\epsilon$ such that
    \[
    G + \epsilon H_\chi \ge 2Ar \cdot \rho^*\alpha(\partial_t) + Ar^2 \cdot \rho^*d\alpha(\partial_t,\partial_r) \ge \epsilon Ar \qquad\text{for all $r < \delta$}
    \]
On the other hand, $G$ is uniformly lower bounded outside of $\R/\Z \times \D(\delta)$, so may choose $\epsilon$ small enough so that $G + \epsilon H_\chi > 0$ away from $\R/\Z \times \D(\delta)$ (and thus for all $r > 0$). 
 \end{proof} 

\noindent The following technical lemma will be useful in a later section.

\begin{lemma} \label{lem:volume_inequality} For any smoothing function $\chi$ and any $C$ with $0 < C < 1$, we have
\begin{equation} \label{eq:}
\alpha_{\epsilon\chi} \wedge d\alpha_{\epsilon\chi} \ge C \cdot \rho^*\alpha \wedge \rho^*d\alpha \qquad\text{for all $\epsilon$ sufficiently small.}
\end{equation}
\end{lemma}

\begin{proof} Note that $\alpha_{\epsilon\chi} \wedge d\alpha_{\epsilon\chi} \to \rho^*\alpha \wedge \rho^*d\alpha$ as $\epsilon \to 0$ away from $\Gamma$, so the result holds away from a small neighborhood of $\Gamma$. On the other hand, as in Lemma \ref{lem:the_set_of_smoothing_functions} we have
\[
\alpha_{\epsilon\chi} \wedge d\alpha_{\epsilon\chi} = (G + \epsilon H_\chi) \cdot dt \wedge dr \wedge d\theta  \qquad\text{and}\qquad \rho^*\alpha \wedge \rho^*d\alpha = G \cdot dt \wedge dr \wedge d\theta
\]
Since $H_\chi$ is positive near $0$, we conclude that $\alpha_{\epsilon\chi} \wedge d\alpha_{\epsilon\chi} \ge \rho^*\alpha \wedge d\alpha$ in a neighborhood of $\Gamma$ independent of $\epsilon$. This proves the result. \end{proof}

An analogous argument can be used to establish the following parametric version of Lemma \ref{lem:volume_inequality}.

\begin{lemma} \label{lem:volume_inequality_parametric} Fix two smoothing functions $\chi$ and $\zeta$ and a constant $C$ with $0 < C < 1$. Suppose that
\[
\alpha_\chi \wedge d\alpha_\chi > C \cdot \rho^*\alpha \wedge \rho^*d\alpha
\]
Then for all $\epsilon$ sufficiently small, the 1-parameter family of smoothing functions  $\chi_s = (1 - s)\chi + s\delta \zeta$ satisfies the same volume form inequality
\begin{equation} \label{eq:}
\alpha_{\chi_s} \wedge d\alpha_{\chi_s} > C \cdot \rho^*\alpha \wedge \rho^*d\alpha \qquad\text{for all $\epsilon$ sufficiently small.}
\end{equation}
\end{lemma}

Recall from Lemma \ref{lem:whitney_contact_structure} that any Whitney contact form has an associated smooth contact structure (well-defined up to isotopy) acquired as the kernel of any smooth (contact) 1-form sufficiently close in the Whitney topology on Whitney forms.

\begin{definition} The \emph{associated contact structure} $\xi$ to a contact form $\alpha$ with singularities along a link $\Gamma$ is defined to be the contact structure associated to the Whitney contact form
\[
\alpha_\chi \qquad\text{for any smoothing chart $(U,\iota,g)$ and smoothing function $\chi$}
\]\end{definition}

\begin{lemma} The associated contact structure $\xi$ to a contact form $\alpha$ with singularities along a link $\Gamma$ is independent of the choice of smoothing chart and smoothing function up to isotopy of contact structures. 
\end{lemma}

\begin{proof} Denote the space of tuples $(U,\iota,g,\chi)$ where $(U,\iota,g)$ is a smoothing chart and $\chi$ is a smoothing function by $\mathcal{S}$. Any continuous family of smoothing charts and smoothing functions induces a continuous family of Whitney contact forms, and thus of corresponding smoothed contact structures. Thus the induced contact structure depends only on the component of $\mathcal{S}$ containing $(U,\iota,g,\chi)$. On the other hand, Lemma \ref{lem:the_set_of_smoothing_functions} implies that $\mathcal{S}$ is homotopy equivalent to the space of smoothing charts $(U,\iota,g)$. Moreover, the space of choices of $g$ is contractible, and the space of tubular neighborhoods is homotopy equivalent to the space of framings, i.e. bundle isomorphisms $\tau$ from $\nu\Gamma$ to the trivial bundle $\Gamma \times \R^2$. 

\vspace{3pt}

This shows that the contact structure only depends (up to isotopy) on the isotopy class of framing induced by $\iota$. To show independence of the framing, assume for simplicity that $\Gamma$ is connected and choose a chart
\[
\iota:\R/\Z \times \mathbb{D} \xrightarrow{\sim} U \subset Y \qquad\text{restricting to a diffeomorphism}\qquad \R/\Z \times 0 \xrightarrow{\sim} \Gamma
\]
This induces a framing and any other isotopy class of framing is induced by a chart of the form
\[
\jmath = \iota \circ \phi_k \qquad\text{where}\qquad \phi_k(t,r,\theta) = (t,r,\theta + 2\pi k t)
\]
The contact forms $\alpha_\iota$ and $\alpha_\jmath$ induced by the tuples $(U,\iota,g,\chi)$ and $(U,\jmath,g,\chi)$ are given in the chart $\iota$  on $U$ by
\[
\alpha_\iota = \rho^*\alpha + \chi(r) d\theta \qquad\text{and}\qquad \alpha_\jmath = \rho^*\alpha + \chi(r) d\theta + 2\pi k\chi(r)dt
\]
Consider the family of 1-forms interpolating between $\alpha_\iota$ and $\alpha_\jmath$ given by
\[
\alpha_s = \rho^*\alpha + \chi(r)(d\theta + 2\pi k s dt) = \rho^*\alpha + \chi(r) d\theta_s \qquad\text{where}\qquad d\theta_s = d\theta + 2\pi k s dt
\]
Let $G$ be the function as in Remark \ref{rmk:coordinate_formulas}. Then we may write
\[
\alpha_s \wedge d\alpha_s = (G + \chi'(r)\rho^*\alpha(\partial_t - 2\pi ks \partial_\theta) + \chi(r)\rho^*d\alpha(\partial_t - 2\pi k s\partial_\theta,\partial_r)) \cdot dt \wedge dr \wedge d\theta_s
\]
Note that $\partial_\theta$ vanishes at $0$ and $\alpha(\partial_t) > 0$ when $r = 0$. Moreover, $\rho^*d\alpha$ is bounded since $\alpha$ is Lipschitz. Thus by scaling $\chi$ as in Lemma \ref{lem:the_set_of_smoothing_functions}, we may assume that there is a small $\delta$ and constant $A,C > 0$ such that $\alpha_s \wedge d\alpha_s$ is a volume form for $r \ge \delta$ and such that, for $r \le \delta$, we have
\[
\chi(r) = Ar^2 \qquad \rho^*\alpha(\partial_t - 2\pi ks\partial_\theta) \ge C \qquad |\rho^*d\alpha(\partial_t - 2\pi k s\partial_\theta,\partial_r)| \le C
\]
In particular, this implies that there is a constant $C' > 0$ such that for all $r \le \delta$, we have
\[
G + \chi'(r)\rho^*\alpha(\partial_t - 2\pi ks \partial_\theta) + \chi(r)\rho^*d\alpha(\partial_t - 2\pi k s\partial_\theta,\partial_r) \ge G + C(2Ar - Ar^2) \ge C'r
\]
As in Lemma \ref{lem:the_set_of_smoothing_functions}, this implies that $\alpha_s \wedge d\alpha_s$ is a volume form for $r \le \delta$ (and thus everywhere) for each $s$. Thus the family $\alpha_s$ is contact for all $s$ and the corresponding contact structures are isotopic.\end{proof}

\subsection{Lefschetz Numbers And Tracking} \label{subsec:lefschetz_and_tracking} In this part, we briefly review versions of the Lefschetz number for homeomorphisms and flows. We also prove a certain stability property of this number for flows. This will be useful in our analysis of the Reeb flow of smoothings.

\vspace{3pt}

We start by recalling the classical Lefschetz index for an isolated fixed point. Fix a continuous map of a manifold to itself 
\[
\Phi:X \to X \qquad \text{with an isolated fixed point $x$}
\]

\begin{definition}[Lefschetz Index] \label{thm:Lefschetz_index} The \emph{Lefschetz index} $\on{Lef}(\Phi;x)$ of an isolated fixed point $x$ is defined as follows. Fix a coordinate chart near $x$, so that $\Phi$ may be viewed as a map
\[\Phi:U \to V\]
of open sets in $\R^n$. Then the Lefschetz index is the degree of the map from the sphere of radius $\epsilon$ around $x$ to the unit sphere, defined by taking $y$ to the unit vector parallel to $\Phi(y) - y$.
\end{definition}

The Lefschetz index is independent of the chart and the parameter $\epsilon$. There is also the following global version for the whole map.

\begin{definition}[Lefschetz Number] The \emph{Lefschetz number} $\on{Lef}(\Phi)$ of a continuous map $\Phi$ with isolated fixed points is the sum
\[
\on{Lef}(\Phi) = \sum_{x \in \on{Fix}(\Phi)} \on{Lef}(\Phi,x)
\]
\end{definition}

\begin{theorem}[Lefschetz-Hopf] \label{thm:Lefschetz_Hopf} The Lefschetz number $\on{Lef}(\Phi)$ of a continuous map $\Phi:X \to X$ of a closed manifold $X$ with finite fixed point sets is invariant under homotopy.
\end{theorem}

\noindent We will require the following local and non-compact version of the Lefschetz-Hopf theorem.

\begin{theorem} \label{thm:rel_LefHopf} Let $\Phi,\Psi:\R^n \to \R^n$ be proper continuous maps with finite fixed point sets such that $\Phi = \Psi$ away from a compact set $K$. Then the Lefschetz numbers coincide.
\[\on{Lef}(\Phi) = \on{Lef}(\Psi)\]
\end{theorem}
\begin{proof} Since $\Phi$ and $\Psi$ are proper and continuous, both maps extend to continuous maps $\bar{\Phi}$ and $\bar{\Psi}$ of the one point compactification $S^n = \R^n \cup \infty$ that agree in a neighborhood of $\infty$. The compactified maps have a single additional fixed point at $\infty$ and so the Lefschetz numbers are
\[
\on{Lef}(\bar{\Phi}) = \on{Lef}(\Phi) + \on{Lef}(\bar{\Phi},\infty) \qquad\text{and}\qquad \on{Lef}(\bar{\Psi}) = \on{Lef}(\Psi) + \on{Lef}(\bar{\Psi},\infty) 
\]
Moreover, the homotopy $\Phi_s = (1-s)\Phi + s\Psi$ from $\Phi$ to $\Psi$ extends to a homotopy $\bar{\Phi}_s$ from $\bar{\Phi}$ to $\bar{\Psi}$. Thus the result follows from Theorem \ref{thm:Lefschetz_Hopf}. \end{proof}

Next, we consider a variant of the Lefschetz number for continuous flows on manifolds. Fix a smooth manifold $X$ and a Lipschitz, fixed-point free flow
\[
\Phi:\R \times X \to X \qquad\text{generated by a vector field}\qquad V
\]
We fix some terminology regarding holonomy maps associated to a local section. 

\begin{definition} A \emph{disk-like local section} $Q$ of $\Phi$ is a  codimension one embedded disk $Q \subset X$ transverse to the flow $\Phi$. The \emph{holonomy map} of $Q$ is the (partially defined) continuous map 
\[
\on{Hol}_\Phi:Q \to Q 
\]
This map sends a point $y \in Q$ to $\Phi(T,y) \in Q$ where $T$ is the minimum time $t > 0$ such that $\Phi(t,x) \in Q$, if there is such a finite $T$, and is otherwise not defined. The function sending a point to this minimal time is the \emph{return time} function
\[
\tau_\Phi:Q \to (0,\infty]
\]
Given an open subset $P \subset Q$ on which the holonomy is defined, the \emph{flow tube} $U(\Phi,P)$ associated to $P$ is the open neighborhood of $\Gamma$ given by
\[
U(\Phi,P) = \on{int}\{\Phi(t,x) \; : \; x \in P \text{ and }0 \le t \le \tau_\Phi(x)\} \subset X
\]\end{definition}

\begin{remark} Note that if a local section $Q$ interects a simple closed orbit $\eta$ of $\Phi$ at a single point $x$, then there is a neighborhood $P$ of $x$ in $Q$ such that the holonomy is well-defined on $P$. Moreover, $x$ is an isolated fixed point of $\on{Hol}_\Phi$ if $\eta$ is an isolated closed orbit. \end{remark}

\begin{definition}[Orbit Index] \label{def:orbit_index} The \emph{Lefschetz index} $\on{Lef}(\Phi,\eta)$ of a (possibly multiply covered) isolated periodic orbit $\eta$ that is a $k$-fold cover of a simple orbit $\gamma$ is defined by
\[  
\on{Lef}(\Phi,\eta) = \on{Lef}(\on{Hol}_\Phi^k,x)
\]
where $\on{Hol}_\Phi:Q \to Q$ is the holonomy map of a cross section intersecting $\gamma$ at $x$. \end{definition}

\begin{example}[Non-Degenerate Orbits] \label{ex:non_deg_Lef} Suppose that $\Phi$ is a smooth (or merely differentiable) flow near an simple orbit $\gamma$ and that $\gamma$ is non-degenerate in the sense that the linearized map
\[
T(\on{Hol}_\Phi)_x:T_xP \to T_xP \qquad\text{at the fixed point $x \cap \gamma$}
\]
has no 1-eigenvalues. In this case, the Lefschetz index of $x$ (and thus of $\eta$) is $\pm 1$ where the sign is given by the sign of the determinant
\[
\on{det}(T(\on{Hol}_\Phi)_x - \on{Id})
\]\end{example}

\begin{example}[Reeb Orbits]  \label{ex:Reeb_Lef}Suppose $\Phi$ is a smooth Reeb flow near an isolated simple orbit $\gamma$ and that $\gamma$ is non-degenerate. Then the Lefschetz index is given by
\[
\on{Lef}(\Phi,\gamma) = (-1)^{|\gamma|}
\]
where $|\gamma| = n - 3 + \CZ(\gamma)$ is the mod 2 SFT grading and the ambient contact manifold is dimension $2n-1$. This follows from the standard relation between the sign of $\on{det}(T(\on{Hol}_\Phi)_x - \on{Id})$ and the Conley-Zehnder index (cf. \cite[p.20]{salamon1999lectures} or \cite{g2014}).\end{example}

\begin{example}[Orbits Of Pseudo-Anosovs]  \label{ex:PA_Lef} Suppose that $\Phi$ is a pseudo-Anosov flow and that $\gamma$ is a closed simple closed orbit. Then
\[
\on{Lef}(\Phi,\gamma) = \left\{
\begin{array}{cc}
1 & \text{if $\gamma$ is rotating singular}\\
-p + 1 & \text{if $\gamma$ is non-rotating singular with $n$ prongs}\\
1 & \text{if $\gamma$ is smooth and negative hyperbolic}\\
-1 & \text{if $\gamma$ is smooth and positive hyperbolic}\\
\end{array}\right. 
\]
The two smooth cases follow easily from Example \ref{ex:non_deg_Lef}. The two singular cases follow by computing the Lefschetz index of the unique fixed point of the standard local model of the relevant singular orbit. This computation can be done by finding a smooth map with non-degenerate fixed points
\[
\phi:\R^2 \to \R^2 \qquad\text{agreeing with the standard map $\phi_{n,k,\lambda}$ outside of a compact set}
\]
Such maps are illustrated in e.g. Cotton-Clay \cite{cotton2009symplectic}. In the rotating case, the perturbed map has a single elliptic fixed point. In the non-rotating case, the perturbed map has $-p+1$ positive hyperbolic fixed points. Theorem \ref{thm:rel_LefHopf} then yields the desired formula.\end{example}

Finally, we discuss a certain stability property for the indices of orbits. Consider the case where the periodic orbits of $\Phi$ are isolated. In particular, if $X$ is compact then the set of periodic orbits below any fixed period is finite and the set of all periods
\[
\on{Spec}(\Phi) \subset \R_+
\]
is discrete. Fix pair of flows $\Phi$ and $\Psi$ with isolated periodic orbits on a manifold $X$ generated by Lipschitz vector fields $U$ and $V$ respectively. We adopt the following terminology.

\begin{definition}[Tracking] \label{def:tracking} The flow $\Psi$ \emph{tracks simple orbits of} $\Phi$ \emph{up to period} $L$ if each simple orbit $\eta$ of $\Phi$ with period less than $L$ admits disk-like local sections
\[
P_\eta \subset Q_\eta \qquad\text{with well-defined holonomy maps}\qquad \on{Hol}_\Phi:P_\eta \to Q_\eta
\]
satisfying the following conditions.
\begin{itemize}
\item[(a)] Under an identification $Q_\eta \simeq \mathbb{D}^n$, the holonomy map extends to a continuous map
\[\on{Hol}_\Phi:\R^n \to \R^n \qquad\text{with a unique fixed point}\qquad x_\eta = \eta \cap P_\eta\]
\item[(b)] The flow tubes $U_\eta = U(\Phi,P_\eta)$ are pairwise disjoint.
\item[(c)] The vector fields $U$ and $V$ generating $\Phi$ and $\Psi$ agree outside of the flow tubes $U_\eta$.
\item[(d)] The holonomy $\on{Hol}_\Psi$
is well-defined on $P_\eta$ with return time bounded above by $L$ on $P_\eta$.
\end{itemize}\end{definition}

\begin{lemma} \label{lem:tracking_Lef} Let $\Psi$ be a flow that tracks the simple orbits of $\Phi$ up to period $L$ and fix a primitive homotopy class $\Gamma \in [S^1,X]$. Then
\[
\sum_{\gamma} \on{Lef}(\Phi,\gamma) = \sum_\eta \on{Lef}(\Psi,\eta)
\]
where the left and right sum are taken over all (simple) closed orbits $\gamma$ of $\Phi$ and $\eta$ of $\Psi$ in the homotopy class $\Gamma$ of period $L$ or less, respectively. \end{lemma}

\begin{proof} Fix a closed orbit $\gamma$ of $\Phi$ of period $L$ or less in the homotopy class $\Gamma$. By Definition \ref{def:tracking}(b,c), the holonomy maps $\on{Hol}_\Phi$ and $\on{Hol}_\Psi$ agree near the boundary of each disk $P_\gamma$. Thus after fixing identifications $Q_\gamma \simeq \D^n$ via Definition \ref{def:tracking}(a), we can be extended these holonomy maps to maps
\[
\on{Hol}_\Phi:\R^n \to \R^n \qquad\text{and}\qquad \on{Hol}_\Psi:\R^n \to \R^n
\]
that agree outside of $\D^n$ and that have fixed points contained in $\D^n$. Then by Theorem \ref{thm:rel_LefHopf}
\[
\on{Lef}(\on{Hol}_\Phi,x_\gamma) = \sum_x \on{Lef}(\on{Hol}_\Psi,x)
\]
where the sum on the right is over all fixed points of $\on{Hol}_\Psi$ in $P_\gamma$. By Definition \ref{def:tracking}(b,d), these fixed points are in bijection with the periodic orbits of $\Psi$ in the neighborhood $U_\gamma$ in the homotopy class $\Gamma$, which all have action $L$ or less. Thus by Definition \ref{def:orbit_index}, we have
\[
\on{Lef}(\Phi,\gamma) = \sum_{\eta \subset U_\gamma} \on{Lef}(\Psi)
\]
where the sum on the right is over all closed orbits of $\gamma$ with period less than $L$ contained in $U_\gamma$ in the class $\Gamma$. The lemma follows by taking the sum of the above equality over $\gamma$.\end{proof}

\subsection{Reeb Dynamics of Smoothings} We conclude this section with several results on the dynamics of the Reeb vector field of the smoothed (Whitney) contact forms in Definition \ref{def:smoothing_function}. For the rest of the section, we again fix a closed 3-manifold
\[
Y \qquad\text{with a contact form $\alpha$ with singularities along a link}\qquad \Gamma \subset Y
\] 
We also choose a smoothing chart for $\alpha$, which will remain fixed for the rest of the section.
\[
\iota:\Gamma \times \D \simeq U \subset Y \qquad\text{and}\qquad g:[0,1]_r \to [0,\infty)
\]
Finally, fix a parametrization of each component of $\Gamma$ by $\R/\Z$ so that each component has a neighborhood with coordinates $(t,r,\theta)$.

\vspace{3pt}

We start by computing an explicit formula for the Reeb vector field of the smoothed (Whitney) contact form. We denote this Reeb vector field by
\[
R_\chi \qquad\text{for a given choice of smoothing function}\qquad \chi \in \mathcal{X}
\]
We also denote the corresponding Reeb flow by $\Phi_\chi$. 

\begin{lemma}[Tracking] \label{lem:tracking} Suppose that $\alpha$ has pseudo-Anosov Reeb flow $\Phi$ and fix a period $L$ not in the spectrum $\on{Spec}(\Phi)$. Then there are constants $\epsilon > 0$ such that, for any smoothing function $\chi$ with
\[\on{supp}(\chi) \subset [0,\epsilon]\]
the Reeb flow $\Phi_{\delta\chi}$ tracks the simple orbits of $\Phi$ up to period $L$ for all sufficiently small $\delta > 0$.\end{lemma}

\begin{proof} For each simple orbit $\gamma$ of $\Phi$ of period less than $L$, fix a tubular neighborhood chart
\[
\iota_\gamma:\R/\Z \times \D \to Y \qquad\text{restricting to a parametrization}\qquad \R/\Z \times 0 \to \gamma \subset Y
\]
that agrees with the smoothing chart around each singular orbit. By shrinking these charts, we may assume that $R$ is transverse to the disks $t \times \D$ in each chart and that each chart is disjoint. We break the construction into a number of steps.

\vspace{3pt}

{\bf Step 1.} First, we construct the disk-like sections as in the definition of tracking (Definition \ref{def:tracking}) and verify the required properties Definition \ref{def:tracking}(a-c). We choose disk-like local sections
\[
Q_\gamma = 0 \times \D \qquad\text{and}\qquad P_\gamma = 0 \times \D(\kappa) \subset 0 \times \D \qquad\text{for $\kappa < 1$ small}
\]
For small $\kappa$, the holonomy of $\Phi$ is a well-defined map $P_\gamma \to Q_\gamma$ for the all of the (finitely many) orbits $\gamma$ of period $L$ or less.  Since $\Phi$ is pseudo-Anosov, the holonomy is conjugate to a linear hyperbolic map (if the orbit is non-singular) or the standard pseudo-Anosov map near a singular fixd point (if the orbit is singular). Thus the holonomy $\on{Hol}_\Phi:P_\gamma \to Q_\gamma$ extends to a map
\[
\on{Hol}_\Phi:\R^2 \to \R^2 \qquad\text{with a unique fixed point at $0 = x_\gamma$} 
\]
This verifies Definition \ref{def:tracking}(a). Definition \ref{def:tracking}(b) can be achieved by taking $\delta$ small enough since the flow is continuous. By choosing $\kappa$ and $\epsilon$ small enough, we can also guarantee that
\[
\R/\Z \times \D(\epsilon) \subset U_\gamma \subset \R/\Z \times \D
\]
In particular, this guarantees that $\alpha_{\delta\chi} = \alpha$ and $R_{\delta\chi} = R$ outside of the union of the sets $\R/\Z \times \D(\epsilon)$ for any $\delta$ across all of the tubular neighborhood charts, which yields Definition \ref{def:tracking}(c). 

\vspace{3pt}

{\bf Step 2.} We now show that the holonomy of the flow of $R_{\delta\chi}$ satisfies Definition \ref{def:tracking}(d) when $\delta$ is sufficiently small. We break this into a few smaller steps.

\vspace{3pt}

{\bf Step 2a.} First, we proving a useful technical claim. Specifically, we claim that for any constant $C$ with $0 < C < 1$, there exists a constant $c > 0$ such that
\begin{equation} \label{eq:tracking_1}
\frac{\chi'}{H_\chi} \ge C \cdot \rho^*dt(R) \qquad\text{in the region $\R/\Z \times \D(c)$}
\end{equation}
Here $dt(R)$ is the $t$-coordinate of $R$ and $H_\chi$ is the function appearing in (\ref{eq:volume_form}). To show this claim, first note that the Reeb vector field is parallel to $\Gamma = \R/\Z \times 0$ and is thus a positive multiple of $\partial_t$ when $r = 0$. In particular
\[
\alpha(\partial_t) = \frac{1}{dt(R)} \cdot \alpha(dt(R) \cdot \partial_t) = \frac{1}{dt(R)} \cdot \alpha(R) = \frac{1}{dt(R)} \qquad\text{along $\Gamma = \R/\Z \times 0$}
\]
This implies that for any $C'$ with $C < C' < 1$, we may choose a small radius $c > 0$ such that
\[
\frac{1}{\rho^*(\alpha(\partial_t))} \ge C' \cdot \rho^*dt(R) \qquad\text{on the region}\qquad \R/\Z \times \D(c)
\]
Now we may expand $H_\chi$ using the formula in (\ref{eq:volume_form}) and get
\[
\frac{\chi'(r)}{H_\chi}  = \frac{1}{\rho^*\alpha(\partial_t)} \cdot \Big(\frac{1}{1 + f}\Big) \quad\text{where}\quad f = \frac{\chi'(r)}{\chi(r)} \cdot \rho^*\Big(\frac{d\alpha(\partial_t,\partial_r)}{\alpha(\partial_t)}\Big)
\]
Now we simply note that for $c$ sufficiently small and all radii $r \le c$ we have
\[
\frac{\chi'(r)}{\chi(r)} = \frac{2Ar}{Ar^2} = 2Ar \qquad d\alpha(\partial_t,\partial_r) \text{is bounded above}\quad\text{and}\quad \alpha(\partial_t) \text{ is bounded below} 
\]
It follows that for any $C''$ with $0 < C'' < 1$, we may choose a $c > 0$ such that for all $r \le c$ we have
\[
\frac{\chi'}{H_\chi} \ge \frac{C''}{\rho^*\alpha(\partial_t)} \ge C'C'' \rho^*dt(R)
\]
Taking $C'$ and $C''$ sufficiently close to $1$ then proves (\ref{eq:tracking_1}).

\vspace{3pt}

{\bf Step 2b.} Next, we claim that for any $C$ with $1 > C > 0$, we have
\begin{equation} \label{eq:t_lower_bound}
dt(R_{\delta\chi}) \ge C \cdot \rho^*(dt(R)) \qquad\text{for sufficiently small $\delta$ in the chart $U = \R/\Z \times \D$} 
\end{equation}
To show this, note that $R_{\delta\chi}$ converges in the $C^\infty$-topology to the Reeb vector field $\rho^*R$ of $\rho^*\alpha$ away from the singularity link $\Gamma$. Thus (\ref{eq:t_lower_bound}) holds on the region $\R/\Z \times \D(c) \subset U$ where $r \ge c$ for any fixed $c$ independent of $\delta$ and $\delta$ small. On the other hand, choose $c > 0$ so that (\ref{eq:tracking_1}) holds. Using the formulas in Remark \ref{rmk:coordinate_formulas}, it is simple to check that (\ref{eq:tracking_1}) implies that
\[
dt(R_\chi) = \frac{\rho^*d\alpha(\partial_r,\partial_\theta) + \chi'(r)}{G + H_\chi} \ge  \frac{\rho^*d\alpha(\partial_r,\partial_\theta)}{G}\cdot (1 + (C - 1) \cdot E) = dt(R) \cdot (1 + (C-1) \cdot E)\]
for any $C$ with $0 < C < 1$. Here $E$ is the function given by
\[
E = \frac{H_\chi}{G/\delta + H_\chi}
\]
Since $G$ is positive away from $r = 0$ and $H_\chi$ is positive near $r = 0$, $E$ achieves its maximum in the region where $H_\chi$ is positive for sufficiently small $\delta$. It follows that
\[
E = \frac{H_\chi}{G/\delta + H_\chi} \le \frac{H_\chi}{H_\chi} \le 1
\]
Therefore it follows that $dt(R_\chi) \ge C \cdot dt(R)$ for sufficiently small $\delta$.

\vspace{3pt}

{\bf Step 2c.} Finally, we demonstrate Definition \ref{def:tracking}(d). This property is automatically satisfied near the smooth orbits of period less than $L$. Thus fix a singular orbit $\gamma$ of $\Phi$ of period $L$. By Step 2b, we can choose $\delta$ so that the Reeb vector field $R_{\delta \chi}$ is everywhere positively transverse to the disks $t \times \D \subset \R/\Z \times \D = U$. This implies that every point $x \in P_\gamma$ has a finite return time to $Q_\gamma$, so that the holonomy $\on{Hol}_\Psi$ is well-defined on the disk-like local section $P_\gamma$. Now let
\[
\eta:[0,1]_s \to U
\]
be a parametrization of the trajectory of $R_{\delta \chi}$ starting at $x$ of length $\tau_\Psi(x)$, parametrized so that $dt(\frac{d\eta}{ds}) = 1$. It is simple to see that we may recover $\tau_\Psi(x)$ as the integral
\begin{equation} \label{eq:tracking_2}
\tau_\Psi(x) = \int_0^1 \frac{1}{dt(R_{\delta\chi}) \circ \eta(s)} ds
\end{equation}
On the other hand, by choosing the surface of section $P_\gamma$ to be sufficiently small in Step 1, we can make $\eta$ arbitrarily close in $C^0$ to the parametrization $\R/\Z \to \gamma$ of $\gamma$ with $dt(\frac{d\gamma}{ds}) = 1$. Then using the continuity of $\rho^*dt(R)$, we may guarantee that
\begin{equation} \label{eq:tracking_3}
\Big|\frac{1}{\rho^*dt(R) \circ \eta(s)} - \frac{1}{\rho^*dt(R) \circ \gamma(s)}\Big| \le b \qquad\text{for arbitrarily small $b > 0$}
\end{equation}
Applying the formulas (\ref{eq:tracking_2}-\ref{eq:tracking_3}) and by choosing $\delta$ sufficiently small as in Step 2b, we find that for any $b > 0$ and $C < 1$, we have
\[
\tau_\Psi(x) = \int_0^1 \frac{1}{dt(R_\chi) \circ \eta(s)} ds \le \frac{1}{C} \cdot \int_0^1 \frac{1}{dt(R) \circ \eta(s)} ds \le \frac{1}{C} \cdot \Big(\frac{1}{dt(R) \circ \gamma(s)} ds - b\Big)= \frac{1}{C}(T - b)
\]
Since the period $T$ is less than $L$, we may choose $C$ and $b$ so that $\tau_\Psi < L$. By doing this over all of the (finitely many) singular orbits, we obtain Definition \ref{def:tracking}(d).\end{proof}

\section{Contact Homology} \label{sec:contact_homology}

In this section, we review a version of contact homology graded by free homotopy classes following \cite{pardon2019contact}. For the rest of the section, fix a closed contact manifold
\[
(Y,\xi)
\]

\begin{remark}[Foundations] Classical transversality methods can be used to construct cylindrical contact homology in the hypertight setting, as demonstrated by the work of Hutchings-Nelson \cite{hutchings2017cylindrical,hutchings2022s} or Bao-Honda \cite{bao2018definition}. This can also be done using more advanced transversality techniques (e.g. VFC methods \cite{pardon2019contact}). It will be convenient for us to base our discussion on \cite{pardon2019contact}.\end{remark}

\subsection{Preliminaries} \label{subsec:preliminaries} We start with some preliminary terminology and definitions on Reeb orbits and contact forms following \cite{pardon2019contact}. 

\vspace{3pt}

We first discuss Reeb orbits, their Conley-Zehnder indices and non-degeneracy. Fix a periodic (and possibly multiply covered) Reeb orbit
\[\Gamma \qquad\text{of a contact form $\alpha$ on $(Y,\xi)$}\]
Let $L$ denote the period of the Reeb orbit. Recall that $\Gamma$ is \emph{non-degenerate} if the linearized return map of the Reeb flow $\Phi$ restricted to $\xi$ satisfies
\[
T\Phi_L|_\xi:\xi_P \to \xi_P\qquad \text{has no $1$-eigenvalues for every $P \in \Gamma$}
\]
Every non-degenerate Reeb orbit $\Gamma$ has a well-defined Conley-Zehnder index and SFT grading \cite[Def 2.48]{pardon2019contact} valued in $\Z/2$. We denote these invariants by
\[
\CZ(\Gamma) = \CZ(\Gamma,\tau) \in \Z/2 \qquad\text{and}\qquad |\Gamma| = n - 3 + \CZ(\Gamma) \in \Z/2
\]
Here $\CZ(\Gamma,\tau)$ denote the Conley-Zehnder index of the (lifted) linearized return map in an arbitrary trivialization of $\Gamma^*\xi$ \cite{g2014} and the dimension of $Y$ is $2n-1$. To each basepoint $P \in \Gamma$, one can assign a rank one graded $\Z$-module, the \emph{orientation line}
\[
\mathfrak{o}_{\Gamma,P} \qquad\text{concentrated in grading }|\Gamma|
\]
The orientation lines naturally form a rank one local system on $\Gamma$ and $\Gamma$ is \emph{good} if this local system is trivial \cite[Def 2.49]{pardon2019contact}. Thus any good orbit has a natural orientation line $\mathfrak{o}_\Gamma$ independent of $P$. 

\vspace{3pt}

We will also need a variant of non-degeneracy for contact forms. Fix a free homotopy class $\Gamma \in [S^1,Y]$ and an action $L \in (0,\infty]$.

\begin{definition} \label{def:non_deg_up_to_L} A contact form $\alpha$ is \emph{non-degenerate in class $\Gamma$ up to action $L$} if every closed orbit of period less than $T < L$ that is either contractible or in the homotopy class $\Gamma$ is non-degenerate.
\end{definition}

Finally, we require a certain directed category (and in fact, a poset) of contact forms. Recall that there is a natural partial order on contact forms for $\xi$ defined by
\[
\alpha > \beta \qquad\text{if}\qquad \alpha = f\beta \text{ for a smooth function $f$ with $f > 1$.}
\]

\begin{definition}[Category Of Pairs] \label{def:cat_of_pairs}  Let $\mathcal{C}(Y,\xi)$ be the category of pairs $(\alpha,L)$ where $\alpha$ is a contact form on $(Y,\xi)$ and $L$ is a positive real number. We define the morphisms by declaring that
\[
(\alpha,L) \to (\beta,M) \qquad\text{exists if}\qquad \alpha/L > \beta/M
\]
and that this morphism is unique. This category is filtered (and in fact a poset).\end{definition}

\noindent We will frequently use this category to define variants of contact homology by taking a colimit with respect to a functor from this category. In particular, the following variant of cofinality will be useful.

\begin{definition}[Strongly Cofinal] \label{def:strongly_cofinal} A sequence of pairs $(\alpha_i,L_i)$ in $\mathcal{C}(Y,\xi)$ is called \emph{strongly cofinal} if there are constants $D > C > 1$ and a fixed contact form $\alpha$ such that
\[
C^{-1} \alpha < \alpha_i < C\alpha \qquad\text{and}\qquad L_{i+1} > D^2L_i
\]
\end{definition}

\begin{lemma}[Cofinality] \label{lem:cofinality} An strongly cofinal sequence $(\alpha_i,L_i)$ is cofinal in $\mathcal{C}(Y,\xi)$.
\end{lemma}

\begin{proof} It follows from the definition that $\alpha_{i+1} < C^2\alpha_i$. Te
\[
\frac{1}{L_i}\alpha_i > \frac{1}{C^2L_i}\alpha_{i+1} > \frac{D}{C^2L_{i+1}}\alpha_{i+1} \qquad\text{and thus}\qquad  \alpha_1 > L_1(D/C^2)^i \cdot \frac{\alpha_i}{L_i}
\]
Since $D > C^2 > 1$, this implies that there is a morphism $(\alpha_i,L_i) \to (\alpha_{i+1},L_{i+1})$ and for any pair $(\beta,K)$, there is a morphism $(\beta,K) \to (\alpha_i,L_i)$ for large enough $i$. \end{proof}

\subsection{Contact Homology Chain Complexes} \label{subsec:CH_in_a_class} We next discuss the chain level constructions for the variants of contact homology that we will need in this paper. Fix a homotopy class
\[
\Gamma \in [S^1,Y]
\]
\begin{definition}[Contact Homology Data] A choice of \emph{contact homology data} for the contact manifold $(Y,\xi)$ in the homotopy class $\Gamma$ is a tuple
\[(L,\alpha,J,\theta)\]
consisting of a choice of action $L \in (0,\infty]$, a contact form $\alpha$  that is non-degenerate in class $\Gamma$ up to action $L$, a translation invariant SFT-type almost complex structure
\[
J \qquad\text{on the symplectization}\qquad \R \times Y
\]
as in \cite[p. 3-4]{pardon2019contact} and virtual fundamental cycle data $\theta$ for the compatified moduli spaces of genus zero SFT-type $J$-holomorphic curves in $\R \times Y$ with one positive output and many negative outputs, where every output and input is in the homotopy class $\Gamma$ (see \cite[Thm 1.1]{pardon2019contact}). \end{definition}

\begin{construction}[Chain Complexes] \label{con:CH_chain_complex} The \emph{contact homology chain complex}
\[
C_L(\alpha,J,\theta;\Gamma) \qquad\text{with differential}\qquad \partial:C_L(\alpha,J,\theta;\Gamma)  \to C_L(\alpha,J,\theta;\Gamma) 
\]
associated to a choice of contact homology data $(L,\alpha,J,\theta)$ is constructed as follows. Consider the graded $\Q$-modules
\[
U_L(\alpha,J,\theta)\qquad\text{and}\qquad V_L(\alpha,J,\theta;\Gamma) 
\]
that are freely generated by the orientation lines of the contractible good Reeb orbits of period less than $L$, and freely generated by the orientation lines of the Reeb orbits of period less than $L$ in the homotopy class $\Gamma$, respectively. Denote the free graded-symmetric algebra generated by the module of contractible orbits by
\[
A_L(\alpha,J,\theta) = SU_L(\alpha,J,\theta)
\]
The contact homology chain groups are given by the following graded tensor product over $\Q$.
\[
C_L(\alpha,J,\theta) = A_L(\alpha,J,\theta) \otimes V_L(\alpha,J,\theta;\Gamma)
\]

The differential $\partial$ is constructed as follows. Given a closed Reeb orbit $\gamma$ and a multiset of good Reeb orbits $\Xi = (\eta_1 \dots \eta_k)$ of the contact form $\alpha$ with period less than $L$, let
\[
\bar{\mathcal{M}}(\gamma,\Xi;A)/\R
\]
denote the compactified moduli space, genus zero, SFT-type pseudo-holomorphic curves in the symplectization $\R \times Y$ with respect to the almost complex structure $J$ modulo translation, in a particular homotopy class $A$ of homology class in $H_2(Y;\gamma \cup \Xi)$. There is an associated Fredholm index (or virtual dimension)
\[
\on{ind}(\gamma,\Xi;A) \in \Z \qquad\text{with}\qquad \on{ind}(\gamma,\Xi;A) = |\gamma| - |\Gamma| \mod 2
\]
As described in \cite[\S 1.2]{pardon2019contact} and \cite[Thm 1.1]{pardon2019contact}, the choice of VFC data $\theta$ determines an element
\[
\#_\theta(\bar{\mathcal{M}}(\gamma,\Xi;A)/\R):\mathfrak{o}_\gamma \otimes \Q \to \mathfrak{o}_\Xi \otimes \Q \qquad\text{in the case where}\qquad \on{ind}(\gamma,\Xi;A) = 1
\]
where $\mathfrak{o}_\Xi$ is the tensor product of the orientation lines of $\gamma_1 \dots \gamma_k$. This can be interpreted as a \emph{virtual count} of points in the compactified moduli space. For a contractible good Reeb orbit $\gamma$ of period less than $L$, we then define
\begin{equation} \label{eq:CH_differential}
\partial|_{\mathfrak{o}_\gamma} = \sum_{\Xi,A} \#_\theta(\bar{\mathcal{M}}(\gamma,\Xi;A)/\R)
\end{equation}
where the sum is over all orbit multisets $\Xi$ of contractible good Reeb orbits and all homology classes $A$ with $\on{ind}(\gamma,\Xi;A) = 1$. Similarly, if $\gamma$ is a closed orbit in the homotopy class $\Gamma$, we define the differential by the same formula where the sum is taken over all orbit multisets $\Xi$ where one orbit is in the homotopy class $\Gamma$ and all of the other orbits are contractible. By linearity, this defines  a map
\[
U_L(\alpha,J,\theta) \oplus V_L(\alpha,J,\theta;\Gamma) \to A_L(\alpha,J,\theta) \otimes V_L(\alpha,J,\theta;\Gamma) = C_L(\alpha,J,\theta;\Gamma)
\]
which then extends to a full map on $C_L(\alpha,J,\theta;\Gamma)$ by the Leibniz rule. Note that an analogous construction yields a dg-algebra differential on the algebra $A_L(\alpha,J,\theta)$. 
\end{construction}

\begin{remark} Often, the differential (\ref{eq:CH_differential}) is written with combinatorial factors accounting for e.g. labeling orders of the punctures of the negative punctures. Here we adopt the convention of simply absorbing these terms into the moduli count $\#(\bar{\mathcal{M}}(\gamma,\Xi)/\R)$ to simplify notation. 
\end{remark}

We now review the basic properties of the contact chain groups in Construction \ref{con:CH_chain_complex}. All of these results may be obtained from the constructions of virtual fundamental cycles given by Pardon \cite{pardon2019contact}. First, we note that the complex is well-defined.

\begin{remark}[Rescaling/Inclusion] Note that there is a canonical chain isomorphism
\[
C_L(\alpha,J,\theta;\Gamma) = C_{CL}(C\alpha,J,\theta;\Gamma) \qquad\text{for any constant}\qquad C > 0
\]
induced by the bijection between the Reeb orbits of $\alpha$ and $C\alpha$. There is also a chain inclusion
\[
\iota_L^K:C_K(\alpha,J,\theta;\Gamma) \to C_L(\alpha,J,\theta;\Gamma) \qquad\text{for any $K \le L$}
\] \end{remark}

\begin{theorem} \cite{pardon2019contact} \label{thm:CH_defined} The contractible contact dg-algebra and contact homology complex in Construction \ref{con:CH_chain_complex} is a chain complex. 
\end{theorem}

\begin{proof} As with standard contact homology, this follows immediately from the construction of virtual point counts by Pardon \cite{pardon2019contact} and the corresponding master equation \cite[Thm 1.1(v)]{pardon2019contact}. 

\vspace{3pt}

Specifically, the master equation states that for any $\gamma,\Xi$ and $A$ with $\on{ind}(\gamma,\Xi;A) = 2$, we have
\begin{equation} \label{eq:master_eq}
0 = \sum \#_\theta(\bar{\mathcal{M}}(\gamma,\eta \cup \Xi_1;A_1)/\R) \circ \#_\theta(\bar{\mathcal{M}}(\eta,\Xi_2;A_2)/\R)
\end{equation}
where the sum is taken over all orbits $\eta$, orbit sets $\Xi_1,\Xi_2$ and homology classes $A_1,A_2$ such that $\Xi = \Xi_1 \cup \Xi_2$, $A$ is the concatenation of $A_1$ and $A_2$ \cite[\S 2.1]{pardon2019contact} and $\on{ind}(\gamma,\eta \cup \Xi_1;A_1) = \on{ind}(\eta,\Xi_2;A_1) = 1$. Now we simply note that if $\gamma$ is contractible, and all of the orbits in $\Xi$ are contractible, and
\[
\#_\theta(\bar{\mathcal{M}}(\gamma,\eta \cup \Xi_1;A_1)/\R) \qquad\text{and}\qquad \#_\theta(\bar{\mathcal{M}}(\eta,\Xi_2;A_2)/\R) \qquad\text{are non-zero}
\]
then $\eta$ must be contractible. Indeed, if the above virtual counts are non-zero then each compactified moduli space is non-empty \cite[Thm 1.1(iv)]{pardon2019contact}. The resulting building in $\bar{\mathcal{M}}(\gamma,\eta \cup \Xi_1;A_1)$ determines (up to homotopy) a map from a compact genus zero surface $\Sigma$ to $Y$ with boundary on $\gamma \cup \Xi_1 \cup \eta$. If $\gamma$ and $\Xi_1 \subset \Xi$ have contractible components, then we can cap off $\Sigma$ at the corresponding boundary components to get a disk bounding $\eta$. Thus $\eta$ is contractible. This implies that the master equation (\ref{eq:master_eq}) holds where the sum is only over contractible orbits $\eta$. This is equivalent to the identity
\[
\partial^2 = 0 \qquad\text{as a map}\qquad A_L(\alpha,J,\theta) \to A_L(\alpha,J,\theta) 
\]
A similar argument shows that the differential squares to zero for $C_L(\alpha,J,\theta;\Gamma)$.\end{proof}

Next, note that cobordism maps can be constructed between the contact homology chain groups for different choices of data. Specifically, there is a commutative diagram of maps
\begin{equation} \label{eq:cobordism_map_1}
\begin{tikzcd}
A_K(\alpha,I,\theta) \arrow[d,"\subset"] \arrow[r] & A_L(\beta,J,\varphi) \arrow[d,"\subset"]\\
C_K(\alpha,I,\theta;\Gamma) \arrow[r] & C_L(\beta,J,\varphi;\Gamma)
\end{tikzcd}
\end{equation}
for any choice of contact homology data $(K,\alpha,I,\theta)$ and $(L,\beta,J,\varphi)$ that satisfy
\[
\alpha/K > \beta/L \qquad\text{or equivalently}\qquad (\alpha,K) \to (\beta,L) \text{ in }\mathcal{C}(Y,\xi)
\]
The map (\ref{eq:cobordism_map_1}) can also be constructed when $\alpha = \beta$ and $K \le L$. Analogous to the differential, these chain maps are constructed by a virtual point count of points in the compactified moduli space of SFT-type pseudo-holomorphic curves in the manifold $\R \times Y$ equipped with an SFT-type compatible almost complex structure that interpolates between the almost complex structures $I$ and $J$ at the ends \cite[\S 1.3]{pardon2019contact}. These maps depend on analogous choices of cobordism map data (e.g. almost complex structure and VFC data).

\begin{theorem}[Independence] \label{thm:CCH_indepence} The cobordism maps  are chain maps and the induced maps on homology are independent of choices of cobordism map data. \end{theorem}

\begin{theorem}[Naturality] \label{thm:CCH_naturality} The cobordism chain maps are natural in the sense that the following diagram commutes for all choices of contact homology data.
\[
\begin{tikzcd}
HA_{L_0}(\alpha_0,I_0,\theta_0) \arrow[r,"(\ref{eq:cobordism_map_1})"] \arrow[d] \arrow[rr, bend left=16, "(\ref{eq:cobordism_map_1})"'] & HA_{L_1}(\alpha_1,I_1,\theta_1) \arrow[d] \arrow[r,"(\ref{eq:cobordism_map_1})"] & HA_{L_2}(\alpha_2,I_2,\theta_2) \arrow[d] \\
HC_{L_0}(\alpha_0,I_0,\theta_0) \arrow[r,"(\ref{eq:cobordism_map_1})"] \arrow[rr, bend left=-12, "(\ref{eq:cobordism_map_1})"'] & HC_{L_1}(\alpha_1,I_1,\theta_1) \arrow[r,"(\ref{eq:cobordism_map_1})"] & HC_{L_2}(\alpha_2,I_2,\theta_2) 
\end{tikzcd}
\]
Moreover, the cobordism maps  $HA_K(\alpha,I,\theta) \to HA_L(\alpha,I,\theta)$ and $HC_K(\alpha,I,\theta) \to HC_L(\alpha,I,\theta)$ for $K < L$ is induced by the inclusion map $\iota^K_L$. \end{theorem}

\noindent These independence theorems are proven by constructing appropriate chain homotopies between the cobordism chain maps (\ref{eq:cobordism_map_1}) corresponding to deformations of almost complex structure and to cobordism composition as in \cite[\S 1.4 and 1.5]{pardon2019contact}. The corresonding chain map and chain homotopy identities can be verified by an identical discussion as in Theorem \ref{thm:CH_defined}.

\vspace{3pt}

\subsection{Contact Homology In A Homotopy Class.} We are now able to fully introduce the variants of contact homology that we will use in the rest of the paper.

\begin{definition}[Contact Homology In A Homotopy Class]\label{def:CH_in_homotopy_class} The \emph{contact homology in the homotopy class $\Gamma$} of a pair $(\alpha,L)$ of contact form and period is the colimit
\[
CH_L(Y,\alpha;\Gamma) = \underset{(\beta,K) \to (\alpha,L)}{\on{colim}}\Big( \underset{(J,\theta)}{\on{colim}}\big(HC_K(\beta,J,\theta;\Gamma)\big)\Big)
\]
where in the outer colimit, the contact form $\beta$ is non-degenerate up to action $K$ in the class $\Gamma$ in the sense of Definition \ref{def:non_deg_up_to_L}. These contact homology groups admit an \emph{action filtration} by subgroups
\[
CH_L^K(Y,\alpha;\Gamma) = \on{im}\big(\iota^K_L:CH_K(Y,\alpha;\Gamma) \to CH_L(Y,\alpha;\Gamma)\big) \subset CH_L(Y,\alpha;\Gamma)
\]
Similarly, the \emph{contact homology in the homotopy class $\Gamma$} of the contact manifold $(Y,\xi)$ is the colimit
\[
CH(Y,\xi) = CH(Y,\xi;\Gamma) = \underset{(\alpha,L)}{\on{colim}}\big(CH_L(Y,\alpha;\Gamma)\big)
\]
These contact homology groups admit an analogous action filtration by subgroups $CH^A(Y,\alpha;\Gamma)$.
\end{definition}

\begin{definition}[Contractible Contact Homology] \label{def:contractible_CH} Analogously, the \emph{contractible contact homology} $CHA_L(Y,\alpha)$ is defined by the analogous colimit of unital algebras
\[
CHA_L(Y,\alpha) = \underset{(\beta,K) \to (\alpha,L)}{\on{colim}}\Big( \underset{(J,\theta)}{\on{colim}}\big(HA_K(\beta,J,\theta;\Gamma)\big)\Big)
\]
where in the outer colimit, the contact form $\beta$ has non-degenerate contractible orbits up to action $K$. Similarly, the \emph{contractible contact homology} of the contact manifold $(Y,\xi)$ is the colimit
\[
CHA(Y,\xi) = CHA(Y,\xi) = \underset{(\alpha,L)}{\on{colim}}\big(CHA_L(Y,\alpha)\big)
\]
These contact homology groups also inherit analogous action filtrations. \end{definition}

\begin{remark} \label{rmk:std_def} By Theorems \ref{thm:CCH_indepence} and \ref{thm:CCH_naturality}, if $\alpha$ is non-degenerate up to action $L$ in the class $\Gamma$ in the sense of Definition \ref{def:non_deg_up_to_L}, then the groups in Theorems \ref{def:contractible_CH} and \ref{def:CH_in_homotopy_class} recover the homologies
\[HC_L(\alpha,J,\theta;\Gamma) \qquad\text{and}\qquad HA_L(\alpha,J,\theta) \qquad\text{for any $(J,\theta)$}\]
This is the more standard definition given for many flavors of contact homology (e.g. in \cite{bao2018definition,hutchings2014cylindrical}).  \end{remark}

There are a number of useful elementary properties of these variants of contact homology that we will need later in the paper. First, we have the following algebraic structure.

\begin{lemma}[Module] The contact homology $CH(Y,\xi;\Gamma)$ in the homotopy class $\Gamma$ is a module over the contractible contact homology algebra $CHA(Y,\xi)$ of $(Y,\xi)$. 
\end{lemma}

\noindent Next, recall that a contact manifold is called \emph{algebraically overtwisted} if the full contact homology, given by the homology of the contact dg-algebra generated by all good Reeb orbits \cite{pardon2019contact} vanishes.
\[
HA(Y,\xi) = 0
\]
This overtwistedness criterion is equivalent to a criterion using the contractible contact homology.

\begin{lemma}[Algebraic Overtwistedness] \label{lem:alg_overtwisted} The contractible contact homology $CHA(Y,\xi)$ vanishes if and only if the full contact homology $HA(Y,\xi)$ of $(Y,\xi)$ vanishes.
\end{lemma}

\begin{proof} Fix a non-degenerate contact form $\alpha$, along with an SFT-type almost complex structure $J$ and VFC data $\theta$. There is a natural inclusion of unital dg-algebras
\[
A(Y,\alpha;J,\theta) \to A_{\on{full}}(Y,\alpha;J,\theta)
\]
from the algebra generated by contractible good orbits. The contact manifold $(Y,\xi)$ is algebraically overtwisted if and only if
\[1 \in  A_{\on{full}}(Y,\alpha;J,\theta)\]
is exact. In particular, $1 = \partial x$ where $x$ may be assumed to be an element of the direct sum of the orientation lines $\mathfrak{o}_\gamma \otimes \mathbb{Q}$ where $\gamma$ is a good closed Reeb orbit. We may write $x = x_\circ + y$ where $x_\circ$ and $y$ are respectively elements of the direct sum of orientation lines of contractible and non-contractible orbits. Then the differential $\partial y$ must have no constant terms, since no non-contractible orbit can bound a genus zero SFT building with no negative ends. Thus
\[
\partial x = \partial x_\circ = 1
\]
We have thus shown that $1$ is exact in $A(Y,\alpha;J,\theta)$ if and only if $1$ is exact in $A_{\on{full}}(Y,\alpha;J,\theta)$. Taking homology yields the desired result.\end{proof}

\begin{corollary} Let $(Y,\xi)$ be an algebraically overtwisted contact manifold. Then $CH(Y,\xi;\Gamma) = 0$.
\end{corollary}

\noindent Next, we have the following standard property that states that the contact homology detects orbits in a given homotopy class.

\begin{lemma}[Orbit Existence] \label{lem:orbit_existence} Let $(Y,\xi)$ be a closed contact manifold with contact form $\alpha$ such that
\[
CH_L(Y,\alpha;\Gamma) \neq 0
\]
Then there is a closed Reeb orbit of $\alpha$ in the homotopy class $\Gamma$ of period $L$ or less. In particular, every contact form has a closed orbit in the class $\Gamma$ is
\[
CH(Y,\xi;\Gamma) \neq 0
\]\end{lemma}

\begin{proof} First, suppose that the contact form $\alpha$ is non-degenerate in the homotopy class $\Gamma$ up to action $L$. In that case, if there is no closed orbit in the class $\Gamma$ then
\[
C_L(\alpha,J,\theta;\Gamma) = A_L(\alpha,J,\theta) \otimes V_L(\alpha,J,\theta;\Gamma) = 0 
\]
Thus the result follows from Remark \ref{rmk:std_def}. In general, we can take a cofinal sequence $(\alpha_i,L_i)$ of such pairs converging to a given $(\alpha,L)$ (in the $C^\infty$-topology on contact forms) and apply Arzel\'{a}-Ascoli to the corresponding sequence of orbits to acquire the desired orbit of $\alpha$.\end{proof}

We conclude this part with a short technical lemma that provides a non-vanishing criterion for the contact homology via strongly cofinal sequences.

\begin{lemma}[Non-Vanishing] \label{lem:non_vanishing_CCH} Fix a contact manifold $(Y,\xi)$, a free homotopy class $\Gamma \in [S^1,Y]$ and a period $A > 0$. Let $(\alpha_i,L_i) \in \mathcal{C}(Y,\xi)$ be a strongly cofinal sequence and suppose that
\[
CH_S(Y,\alpha_i;\Gamma) \to CH_{L_i}(Y,\alpha_i;\Gamma)  \text{ is injective and non-zero for all $A \le S < L_i$ and $i$ sufficiently large}
\]
Then $CH(Y,\xi;\Gamma)$ is non-zero. Moreover, for any contact form $\alpha$ with $\alpha < \alpha_i$ for all $i$, we have
\[CH^A(Y,\alpha;\Gamma) \neq 0\]\end{lemma}

\begin{proof} There exists a $C > 0$ such that $\alpha_i < C^2\alpha_j$ for any $i$ and $j$. Thus, for sufficiently large $i$ and sufficiently large $S,T$ independent of $i$ with $A < S < T$, we have a commutative diagram
\[
\begin{tikzcd}
CH_A(Y,\alpha_i;\Gamma) \arrow[r,"(1)"] \arrow[d,"="] & CH_T(Y,\alpha_i;\Gamma) \arrow[r] & CH_{L_i}(Y,\alpha_i;\Gamma)  \arrow[d]\\
CH_A(Y,\alpha_i;\Gamma) \arrow[r,"(2)"] & CH_S(Y,\alpha_j;\Gamma) \arrow[r,"(3)"] \arrow[u] & CH_{L_j}(Y,\alpha_j;\Gamma)
\end{tikzcd}
\]
The composition of the top row maps is injective, so the map (1) is injective. The map (2) is injective since the diagram commutes. Finally, (3) is injective by hypothesis. Thus the composition of the bottom row is injective. By taking a colimit, this yields an injection
\[ 
CH_A(Y,\alpha_i;\Gamma) \to \underset{j}{\on{colim}}\big(CH_{L_j}(Y,\alpha_j;\Gamma)\big) = CH(Y,\xi;\Gamma)
\]
The final claim of the lemma follows from the fact that we can factor the above maps as
\[
CH_A(Y,\alpha_i;\Gamma) \to CH^A(Y,\alpha;\Gamma) \subset CH(Y,\xi;\Gamma) \qedhere
\]\end{proof}

\begin{remark} \label{rmk:grading_nonvanishing} In Lemma \ref{lem:non_vanishing_CCH}, it suffices to assume that the maps $CH_S(Y,\alpha_i;\Gamma) \to CH_{L_i}(Y,\alpha_i;\Gamma)$ are injective in a single grading that is independent of $S$ and $i$.
\end{remark}

\subsection{Cylindrical Contact Homology.}  \label{subsec:asymptotic_hypertightness} There is a special class of contact manifolds for which the contact homology in a homotopy class is particularly simple.

\begin{definition}[Asymptotically Hypertight] \label{def:asymptotically_hypertight}  A pair $(\alpha,L)$ of a contact form $\alpha$ and a period $L$ is \emph{hypertight} if there are no contractible Reeb orbits of the Reeb vector field with action less than $L$. 

\vspace{3pt} 

A contact structure $\xi$ on $Y$ is \emph{asymptotically hypertight} if the inclusion of the full sub-category of hypertight pairs into the category of all pairs
\[
\mathcal{C}_{\on{HT}}(Y,\xi) \subset \mathcal{C}(Y,\xi) \qquad\text{is cofinal.}
\]
\end{definition}

Asymptotic hypertightness implies a strong version of tightness. In particular, we have the following well-known lemma which is folklore. 

\begin{lemma}[Universally Tight] \label{lem:universally_tight} Any asymptotically hypertight contact structure $\xi$ on a closed manifold $Y$ is tight (and universally tight if $Y$ is 3-dimensional).
\end{lemma}

\begin{proof} In this case, there is a cofinal sequence of pairs $(\alpha_i,L_i)$ such that $\alpha_i$ has no contractible orbits of period $L_i$ of less. Therefore
\[
A_{L_i}(Y,\alpha_i) = \mathbb{Q} \qquad\text{and thus}\qquad CHA_{L_i}(Y,\alpha_i) = \mathbb{Q}
\]
Taking the colimit for $i \to \infty$ then shows that $CHA(Y,\xi) \neq 0$. By Lemma \ref{lem:alg_overtwisted}, this implies that $(Y,\xi)$ is algebraically tight (i.e. not algebraically overtwisted). It is well known (cf. \cite{casals2019geometric}) that any overtwisted contact manifold must be algebraically overtwisted, so $(Y,\xi)$ is tight. Note that any finite cover of a asymptotically hypertight contact manifold is also asymptotically hypertight, so this implies that
\[
(P,\pi^*\xi) \text{ is tight for any finite cover $\pi:P \to Y$}
\]
Finally, suppose that $Y$ were 3-dimensional and not universally tight. Then there is an overtwisted disk $D \subset \widetilde{Y}$ in the universal cover. Since $D$ is compact, the set
\[
S = \{g \in \pi_1(Y) \; : \; gD \cap D \neq \emptyset\} \subset \pi_1(Y)
\]
is finite. On the other hand, the fundamental group of any compact 3-manifold is residually finite. This was proven by Hempel \cite{hempel1987residual} for Haken manifolds and the proof extends to all compact 3-manifolds via the Geometrization Theorem. This implies that for any finite set $S \subset \pi_1(Y)$, there is a finite index, normal subgroup $N \subset \pi_1(Y)$ with $N \cap S = 1$. The quotient
\[
P = \widetilde{Y}/N \to Y
\]
is then a closed cover of closed cover. The projection of $D$ to $P$ is then an embedded overtwisted disk in $P$. This contradicts the tightness of all finite covers of $Y$. \end{proof}

\begin{remark}[Cylindrical CH] For an asymptotically hypertight contact manifold $(Y,\xi)$, the contact homology group in a free homotopy class $\Gamma$ is precisely the cylindrical contact homology over $\Q$ in the homotopy class $\Gamma$. Indeed, for any hypertight pair $(\alpha,L)$, there are cylindrical contact homology groups
\[
CCH_L(Y,\alpha;\Gamma) = CH_L(Y,\alpha;\Gamma) \qquad\text{in the homotopy class $\Gamma \in [S^1,Y]$}
\]
computed from a chain complex freely generated by closed orbits in the homotopy class $\Gamma$. The contact homology in the homotopy class $\Gamma$ is then the colimit
\[
CCH(Y,\xi;\Gamma) = \on{colim}_{(\alpha,L)} \; CCH_L(Y,\alpha;\Gamma)
\]
taken over all hypertight pairs with contact form non-degenerate in the class $\Gamma$.\end{remark}

\subsection{Giroux Torsion And CH Bound} \label{subsec:torsion_and_CH} We next discuss a result on the contact homology and Giroux torsion, generalizing a computation due to Bourgeois-Giroux \cite{bourgeois2005homologie} in the hypertight case. We start by recalling the precise definition of Giroux torsion. We require the following notation.

\begin{notation} We denote the standard 2-dimensional torus
\[
\mathbb{T}^2 = \R/\Z \times \R/\Z \qquad\text{with $(x,y)$-coordinates}
\]
The set of free homotopy classes $[S^1,\mathbb{T}^2]$ is identified with $\Z^2$ in the standard way. We denote the projection map from $[S^1,\mathbb{T}^2]$ to the first $\Z$ factor by
\[
\theta:[S^1,\mathbb{T}^2] \to \mathbb{Z} 
\] \end{notation}

\begin{definition}[Giroux Torsion] \label{def:torsion} A contact 3-manifold $(Y,\xi)$ has \emph{Giroux torsion $k$} in the homotopy class of map $T:\mathbb{T}^2 \to Y$ it admits a contact embedding from the \emph{Giroux torsion domain}
\[
U_k = [0,k]_t \times \mathbb{T}^2 \qquad\text{with contact structure}\qquad \alpha_k = \cos(2\pi t)dx + \sin(2\pi t) dy
\]
such that the restriction of $T$ to the torus $0 \times \mathbb{T}^2$ is homotopic to $T$. \end{definition}

The main result of this part can now be stated as follows.

\begin{theorem}[CH And Giroux Torsion] \label{thm:giroux_torsion_CH} Let $(Y,\xi)$ be a closed contact 3-manifold with Giroux torsion $k$ in the class of essential torus embedding $T:\mathbb{T}^2 \to Y$. Fix a primitive free homotopy class
\[
\Gamma = [S^1,\mathbb{T}^2] \qquad\text{such that}\qquad \theta(\Gamma) \le 0
\]
Finally, suppose that $(Y,\xi)$ is algebraically tight. Then the contact homology in the class $T_*\Gamma$ satisfies
\[
\on{rank}(CH(Y,\xi;T_*\Gamma)) \ge 2k
\]
\end{theorem}

\begin{remark} In fact, one can prove the stronger fact that the rank of $CH(Y,\xi;T_*\Gamma)$ as a module over $CHA(Y,\xi)$ is bounded below by $2k$. For simplicity, we will restrict to the above statement.
\end{remark}

\noindent The proof is similar to the hypertight case \cite{bourgeois2005homologie} and the computation of the cylindrical contact homology of the Giroux torsion $k$ tight contact structure on $\mathbb{T}^3$ (cf. Wendl \cite[Ch 10]{w2016}).

\vspace{3pt}

We now proceed to the proof of Theorem \ref{thm:giroux_torsion_CH}. For simplicity, we prove the case of Giroux torsion one, since the proof in the higher torsion case is similar. We will also assume for simplicity that $\theta(\Gamma) < 0$ and remark on the case of $\theta(\Gamma) = 0$ below. Fix a contact 3-manifold
\[
(Y,\xi) \qquad\text{with embedded Giroux domain}\qquad U = U_1 \subset Y
\]
We start by fixing a good choice of contact form and almost complex structure. Let $\beta$ be a contact form for $\xi$ such that
\[
\beta = \alpha_1 \qquad\text{in a neighborhood of $U$} 
\]
Note that for such a contact form, the boundary $\partial U$ is invariant under the Reeb flow. The closed orbits of the contact form $\beta$ in the free homotopy class $\Gamma$ that are contained in $U_k$ form a Morse-Bott torus of closed embedded Reeb orbits
\[
T_\Gamma \subset \on{int}(U)
\]
Finally, every contractible closed orbit is disjoint from $U$. By taking a $C^\infty$-small perturbation of $\beta$ supported away from $U$, we may further assume that the contractible closed orbits and the closed orbits in the homotopy class $\Gamma$ outside of $U$ are non-degenerate. 

\vspace{3pt}

Next, we can perform a standard Morse perturbation of the Morse-Bott contact form $\beta$. This is carried out carefully in Wendl \cite[\S 10]{w2016} (and also see Bourgeois \cite{bourgeois2002morse} and Yao \cite{yao2022cascades}). Precisely, by \cite[Lecture 10, \S 3]{w2016} there is a contact form $\beta$ that is $C^\infty$-close to $\alpha$ such that
\[
\alpha = \beta \qquad\text{outside of a neighborhood of $T_\Gamma$}
\]
Morover, the closed Reeb orbits of $\alpha$ within $U$ in the homotopy class $\Gamma$ are in bijection with the critical points of a perfect Morse function on the quotient $T_\Gamma/\R \simeq S^1$ of the Morse-Bott family $T_\Gamma$ by the Reeb flow. In particular, there are exactly two closed orbits
\[
\hat{\gamma} \qquad\text{and}\qquad \check{\gamma}
\]
in the class $\Gamma$ contained in $U$, corresponding to the maximum and minimum of the chosen Morse function. Moreover, by \cite[Lemma 10.34]{w2016} there exists a choice of SFT-type almost complex structure $J$ such that the SFT-type $J$-holomorphic cylinders
\[
u:\R \times \R/\Z \to \R \times U \qquad\text{asymptotic to one of $\hat{\gamma}$ or $\check{\gamma}$ at $\pm\infty$}
\]
are transversely cut out and in bijection with gradient flowlines of $f$ connecting $\hat{\gamma}$ and $\check{\gamma}$. We now have the following lemma controlling SFT-type holomorphic curves with respect to $J$.

\begin{lemma}[No Crossing Boundary] \label{lem:nocrossing} Consider a genus zero $J$-holomorphic building
\[
u = (u_1,\dots,u_m) \qquad\text{in the symplectization }\R \times Y
\]
asymptotic to a closed orbit $\gamma$ in the homotopy class $\Gamma$ at $+\infty$ and to closed orbits $\gamma_1 \dots \gamma_k$ at $-\infty$, where $\gamma_1$ is in the homotopy class $\Gamma$ and $\gamma_2 \dots \gamma_k$ are contractible. Then every level is disjoint from $\R \times \partial U$.\end{lemma}

\begin{proof} We argue by contradiction. Consider a connected component of some level $u_i$, denoted by
\[
v:\Sigma \to \R \times Y
\]
Note that $v$ is necessarily a genus zero curve asymptotic to a single orbit $\eta$ at $+\infty$ and a set of orbits $\Xi$ at $-\infty$. We partition $\Xi$ into subsets
\[
\Xi = \Xi_U \sqcup \Xi_Y
\]
where $\Xi_U$ is the set of orbits in $U$ and $\Xi_Y$ is the set of orbits in $Y \setminus U$. Moreover, $\eta$ is either contractible or in the homotopy class $\Gamma$. In the former case, every orbit in $\Xi$ is contractible. In the latter case, there is exactly one orbit in the set $\Xi$ that is in the homotopy class $\Gamma$, and every other orbit is contractible. We assume for the rest of the lemma that we are in the latter case (as the former case is similar) and denote the orbit in the homotopy class $\Gamma$ by $\eta_1$.

\vspace{3pt}

We also let $\chi \subset \Sigma$ denote the inverse image $v^{-1}(\R \times \partial U)$. For simplicity, assume that $v$ is transverse to the hypersurface $\R \times \partial U$, so that the inverse image $\chi$ is a 1-manifold that divides $\Sigma$ into two punctured surfaces with compact boundary
\[
\Sigma_U = v^{-1}(\R \times U) \qquad\text{and}\qquad \Sigma_Y = v^{-1}(\R \times (Y \setminus \on{int}(U)))
\]
Finally, note that by construction $\beta|_{\partial U}$ is simply the closed 1-form $dx$. It follows that
\[
\int_\chi u^*\alpha = [dx] \cdot [\chi] = \theta([\Xi])
\]

Now assume for contradiction that $\chi$ is non-empty. We consider several cases. First, consider the case where both $\eta$ and $\eta_1$ are contained in $U$. In this case, it is simple to see that $\chi$ is null-holomogous (and thus null-homotopic) in $U$. Then by positivity of the action of closed orbits and area of SFT-type $J$-holomorphic curves, we have
\[
0 < \int_{\Sigma_Y} d\alpha = \int_\chi \alpha - \int_{\Xi_Y} \alpha = [dx] \cdot [\chi] -  \int_{\Xi_Y} \alpha = - \int_{\Xi_Y} \alpha < 0 \qquad\text{giving a contradiction.}
\]
A similar argument applies in the parallel case where both $\eta$ and $\eta_1$ are contained in $Y \setminus U$. Second, consider the case where $\eta$ is contained in $U$ and $\eta_1$ is contained in $Y \setminus U$. In this case, it is simple to see that $\Xi$ is in the homology class $\Gamma \in H_1(U;\Z) \simeq \pi_1(U)$ and thus
\[
0 < \int_{\Sigma_Y} d\alpha = \int_\chi \alpha - \int_{\Xi_Y} \alpha < \theta([\chi]) < 0\qquad\text{giving a contradiction.}
\]
A similar argument applies in the parallel case where $\eta$ is contained in $Y \setminus U$ and $\eta_1$ is contained in $U$. This proves the lemma in the case where $u$ is transverse to $\R \times \partial U$.

\vspace{3pt}

Finally, the case where $u$ is not transverse to $\R \times \partial U$ can be handled by applying an analogous argument where $\partial U$ is replaced by the hypersurface
\[
\{\epsilon,1-\epsilon\} \times \mathbb{T}^2 \qquad\text{for generic choice of $\epsilon$ small}
\]
and then taking the limit as $\epsilon \to 0$. \end{proof}

We are now ready to complete the proof of Theorem \ref{thm:giroux_torsion_CH} 

\begin{proof} (Theorem \ref{thm:giroux_torsion_CH}) The contact homology in the homotopy class $T_*\Gamma$ can be computed as the homology of the complex
\[
C = C(\alpha,J,\theta;T_*\Gamma)
\]
as constructed in Section \ref{subsec:CH_in_a_class}. By construction, this complex can be written as a tensor product
\[
C = A \otimes V
\]
where $A = SU(\alpha,J,\theta)$ is the free graded-commutative algebra generated by the orientation lines of closed contractible Reeb orbits and $V = V(\alpha,J,\theta;T_*\Gamma)$ is the free vector space generated by orbits in the homotopy class $\Gamma$. By our construction of $\alpha$, we can split $V$ as a direct sum
\[
V = V_U \oplus V_Y
\]
where $V_U$ and $V_Y$ are spanned by orbits in the homotopy class $\Gamma$ contained in $U$ and the complement of $U$, respectively. This yields a  direct sum splitting of $C$ into sub-modules
\[
C = C_U \oplus C_Y \qquad\text{where}\qquad C_U = A \otimes V_U \qquad\text{and}\qquad C_Y = A \otimes V_Y
\]

Next, note that Lemma \ref{lem:nocrossing} implies that this is a splitting of complexes. Indeed, the differential $\partial$ counts points in the compactified moduli space of pseudo-holomorphic curves with positive end on an orbit in the homotopy class $\Gamma$, a negative end on another orbit in the homotopy class $\Gamma$, and all other negative ends contractible. Such a curve (or more generally, building) cannot cross $\R \times \partial U$ by Step 2. A similar application of Lemma \ref{lem:nocrossing} implies that the identification
\[
C_U = A \otimes V_U
\]
is a tensor product of complexes, since all of the orbits generating $A$ are contractible and thus disjoint from $U$. This implies that
\[
CH(Y,\xi;T_*\Gamma) \simeq (H(A) \otimes H(C_U)) \oplus H(C_Y)
\]
Since $Y$ is algebraically tight, we know that the homology $H(A)$ of the dg-algebra of contractible orbits is non-zero. Therefore
\[
\on{rank}\big(CH(Y,\xi;T_*\Gamma)\big) \ge \on{rank}(H(C_U))
\]
Finally, the construction of the contact form $\alpha$ and almost complex structure $J$ implies that
\[
H(C_U) \simeq H(S^1;\mathbb{Q}) \qquad\text{and thus}\qquad \on{rank}\big(CH(Y,\xi;T_*\Gamma)\big) \ge 2
\]
This proves the inequality if $Y$ has Giroux torsion 1 and the higher torsion case is identical. \end{proof}

\begin{remark} In the case where $\theta(\Gamma) = 0$, the proof can be modified as follows. The contact form $\alpha_U$ on $U = [0,1] \times \mathbb{T}^2$ has two Morse-Bott tori of orbits in the homotopy class $\Gamma$, namely the boundary tori $\partial U$. One can extend the embedding of $U$ into $Y$ to a contact embedding from
\[
(-\epsilon,1+\epsilon) \times \mathbb{T}^2 \subset Y
\]
for a sufficiently small $\epsilon$. The existence of this extension follows by applying a standard neighborhood theorem near $\partial U$. One can then use the same constructions to produce a perturbed contact form $\alpha$ on $U$ with two orbits corresponding to the Morse-Bott torus $0 \times \mathbb{T}^2$. An analogue of Lemma \ref{lem:nocrossing} then holds, where $\partial U$ is replaced the boundary of
\[
[-\epsilon,1-\epsilon] \times \mathbb{T}^2 \subset Y
\]
With these modifications in place, the proof provided above carries through.\end{remark}

\subsection{Entropy And CH Growth Rate} We conclude this section with a lower bound for topological entropy via the homotopical growth rate of contact homology, which generalizes similar estimates for cylindrical contact homology due to Alves \cite{alves2016cylindrical}. 

\vspace{3pt}

The topological entropy of a flow on a smooth manifold is closely related to a different quantity called the homotopical growth rate. We start by recalling the definition of this quantity.

\begin{definition}[Homotopical Growth] The \emph{homotopical growth function} $\on{GF}(\Phi;-)$ of a continuous flow $\Phi$ on a topological space $X$ is the function on $\R$ given by
\[
\on{GF}(\Phi;L) = \#\{\Gamma \in [S^1,X] \; : \; \text{$\Gamma$ is the homotopy class of a closed orbit of period $L$ or less} \}
\]
The \emph{homotopical growth rate} $\on{Gr}(X,\Phi)$ is given by the following formula.
\[
\on{Gr}(X,\Phi)= \limsup_{L}\big(\frac{1}{L} \cdot \on{log}(\on{GF}(\Phi;L)) \big) \in [0,\infty]
\]\end{definition}

\noindent The following result is folklore in the dynamics literature, and a proof can be found in Alves \cite{alves2016cylindrical}.

\begin{theorem} \cite[Thm 1]{alves2016cylindrical} \label{thm:alves_entropy} Let $\Phi$ be a differentiable flow on a manifold $X$. Then the topological entropy $\on{Ent}(X,\Phi)$ is bounded below by the homotopical growth rate.
\[
\on{Ent}(X,\Phi) \ge \on{Gr}(X,\Phi)
\]
\end{theorem}

\noindent Likewise, the following result was proven by Barthelm\'{e}-Fenley \cite{barthelme2017counting} in the Anosov case, and the pseudo-Anosov case follows by the same argument (cf. \cite[p. 25]{barthelmé2022orbitequivalencespseudoanosovflows}).

\begin{theorem} \cite[Thm A]{barthelme2017counting} \label{thm:PA_entropy}  Let $\Phi$ be transitive pseudo-Anosov flow on a 3-manifold $Y$. Then the topological entropy $\on{Ent}(Y,\Phi)$ is positive and equal to the homotopical growth rate.
\[
\on{Ent}(Y,\Phi) = \on{Gr}(Y,\Phi) > 0
\]
\end{theorem}

\noindent There is a direct analogue of the homotopical growth rate that is defined using contact homology.

\begin{definition}[CH Homotopical Growth] \label{def:CH_homotopical_growth} The \emph{CH homotopical growth function} $\on{CHF}(\alpha;-)$ of a contact form $\alpha$ on a closed contact manifold $(Y,\xi)$ is the function on $\R$ given by
\[
\on{CHF}(\alpha;L) =  \#\big\{\Gamma \in [S^1,X] \; : \; CH^L(Y,\alpha;\Xi) \neq 0 \text{ for a primitive $\Xi$ with $\Gamma = \Xi^k$}\big\}
\]
Similarly, the \emph{CH homotopical growth rate} $\on{CHGr}(Y,\alpha)$ is given by the following formula.
\[
\on{CHGr}(Y,\alpha) =  \limsup_{L}\big(\frac{1}{L} \cdot \on{log}(\on{CHF}(\alpha;L)) \big) \in [0,\infty]
\]
\end{definition}

We require a few basic properties of the contact homology growth rate. First we have the following properties with respect to the ordering on contact forms.

\begin{lemma} \label{lem:CHG_properties} Fix contact forms $\alpha$ and $\beta$ on $(Y,\xi)$ with $\alpha \le \beta$ and a constant $A > 0$. Then
\[
\on{CHGr}(Y,\alpha) \ge \on{CHGr}(Y,\beta) \qquad \on{CHGr}(Y,A \cdot \alpha) = A^{-1} \cdot \on{CHGr}(Y,\alpha)
\]
\end{lemma}

\begin{proof} These follow easily from the fllowing properties of filtered contact homology.
\[
CH^L(Y,\beta;\Gamma) \subset CH^L(Y,\alpha;\Gamma) \qquad\text{and}\qquad CH^{AL}(Y,A \cdot \alpha;\Gamma) = CH^{L}(Y,\alpha;\Gamma)
\]
In particular, these imply the corresponding properties of the growth functions
\[
\on{CHF}(\alpha;L) \ge \on{CHF}(\beta;L) \qquad\text{and}\qquad \on{CHF}(A \cdot \alpha;L) = \on{CHF}(\alpha;AL)
\]
and the corresponding results for the growth rate follow from the definition.\end{proof}

\noindent We also have the following lemma, immediate from Lemma \ref{lem:orbit_existence} and Theorem \ref{thm:alves_entropy}. 

\begin{lemma} \label{lem:CHG_entropy_bound} Let $(Y,\xi)$ be a closed contact manifold and let $\Phi$ be a Reeb flow for a contact form $\alpha$. Then
\[
\on{Ent}(Y,\Phi) \ge \on{CHGr}(Y,\alpha)
\]
\end{lemma}

\section{Contact Homology Of Pseudo-Anosov Contact Structures} \label{sec:contact_homology_of_PA_structures} 

In this part, we compute the cylindrical contact homology of a pseudo-Anosov contact structure and prove the main results in the introduction.

\subsection{Construction Of A Cofinal Sequence}  Fix a pseudo-Anosov contact 3-manifold
\[
(Y,\xi) \qquad\text{with pseudo-Anosov contact form}\qquad \alpha
\]
Also fix a smoothing chart $(U,\iota,g)$ near the singular orbits (Definition \ref{def:smoothing_chart}). Recall that a smoothing function $\chi$ in the smoothing chart (Definition \ref{def:smoothing_chart}) is a compactly supported smooth function
\[
\chi:[0,1)_r \to \R \qquad \qquad\text{such that $\alpha_\chi$ is contact}
\]
Finally, let $\Phi$ denote the (pseudo-Anosov) Reeb flow of the contact form $\alpha$.

\begin{lemma}[Tracking $\Rightarrow$ Hypertight] \label{lem:tracking_hypertight} Let $\chi$ be a smoothing function such that $\alpha_\chi$ tracks the Reeb flow of $\alpha$ for action less than $L$. Then $\alpha$ is $L$-hypertight.
\end{lemma}

\begin{proof} In this setting, any closed orbit $\eta$ of action bounded by $L$ must either be contained in $\on{supp}(\chi) \subset U$ or is in the complement of $\on{supp}(\chi)$. In the first case, $\eta$ is homotopic to a multiple cover of one of the singular orbits in the link $\Gamma \subset Y$ of singular orbits. In the second case, $\eta$ is a closed orbit of the (unpurturbed) Reeb flow of $\alpha$. In either case
\[
[\eta] \in [S^1,Y]
\]
is the homotopy class of a closed orbit of the pseudo-Anosov Reeb flow $\Phi$ of $Y$ and is thus non-zero by Lemma \ref{lem:non_contractible_orbit}. This shows that $\alpha_\chi$ is $L$-hypertight.\end{proof}

\begin{lemma}[Cofinal For Pseudo-Anosov] \label{lem:cofinal_for_pseudoanosov} Fix any period $A > 0$ not in the spectrum of $\Phi$. Then there exists a sequence of non-degenerate smooth contact forms $\alpha_i$, periods $L_i$ and contactomorphisms
\[\phi_i:(Y,\xi) \to (Y,\on{ker}(\alpha_{\chi_i}))\]
isotopic to the identity, such that the contact forms $\alpha_{\chi_i}$ satisfy
\[
\text{$\alpha_i$ tracks the Reeb flow of $\alpha$ to action $A$ and up to action $L_i$} \qquad\text{and}\qquad \text{$(\phi^*_i\alpha_{\chi_i},L_i)$ \text{ is strongly cofinal}}
\]
 \end{lemma}

\begin{proof} Fix a pair of constants $D > C > 1$ and actions $L_i > A$ not in the spectrum of $\Phi$ such that $L_{i+1} > D^2 L_i$. We assume that $C$ satisfies the following
\[
\on{log}(C) > \underset{U}{\on{max}} \; 4|\rho^*d\theta(R)|
\]By applying Lemma \ref{lem:tracking} twice, we may choose a sequence of smoothing functions $\bar{\chi}_i$ and constants $\epsilon_i$ such that, for any $\delta_i \le \epsilon_i$, the smoothing function $\chi_i = \delta_i \cdot \bar{\chi}_i$ has the property that
\begin{equation} \label{eq:tracking_sequence}\text{$\alpha_{\chi_i}$ tracks the Reeb flow of $\alpha$ to actions $A$ and $L_i$}\end{equation}
We may assume that $|\chi_i| \le 1$ for all $i$ by choosing $\delta_i$ small enough. Moreover, by Lemma \ref{lem:volume_inequality_parametric}, we may choose $\delta_i$ so that the convex combinations
\begin{equation} \label{eq:cofinal_pa_1}
\chi_{i,s} = s\chi_1 + (1 - s)\chi_i \qquad\text{satisfy}\qquad \alpha_{\chi_{i,s}} \wedge d\alpha_{\chi_{i,s}} \ge \frac{1}{2} \cdot \rho^*\alpha \wedge \rho^*d\alpha 
\end{equation}

Note that the contact forms $\alpha_{\chi_i}$ are only Whitney in the sense of Definition \ref{def:whitney} (and not smooth). However, by Lemma \ref{lem:smooth_dense_Whitney} we may approximate these contact forms in the Whitney topology (Definition \ref{def:whitney}) by smooth contact forms $\alpha_i$. Moreover, by a further $C^\infty$-small perturbation, we may assume that the contact forms $\alpha_i$ are non-degenerate. We denote the Reeb vector-fields by $R_i$. Since the Reeb vector-fields $R_i$ approximate the Reeb vector-fields $R_{\chi_i}$ in the $C^0$-topology, we have an analogous tracking property to (\ref{eq:tracking_sequence}) for the contact forms $\alpha_i$.
\[
\text{$\alpha_i$ tracks the Reeb flow of $\alpha$ to actions $A$ and $L_i$}
\]
We may also approximate the family $\alpha_{\chi_{i,s}}$ by a family of smooth contact forms $\alpha_{i,s}$. In particular
\[
R_{i,s} \text{ approximates }R_{\chi_{i,s}} \qquad\text{and}\qquad |\frac{d\alpha_{i,s}}{ds}(R_{i,s})| \le 2|\frac{d\alpha_{\chi_{i,s}}}{ds}(R_{\chi_{i,s}})|
\]

We next construct the desired contactomorphisms $\phi_i$ by utilizing Gray stability. Recall (cf. \cite[\S 2.2]{Geiges_2008}) that Gray stability provides a 1-parameter family of contactomorphisms
\[
\phi_{i,s}:Y \to Y \qquad\text{with}\qquad \phi_{i,s}^*\alpha_{i,s} = F_{i,s} \cdot \alpha_1 \qquad\text{and}\qquad \frac{d}{ds}(\on{log} F_{i,s}) = \frac{d\alpha_{i,s}}{ds}(R_{i,s}) \circ \phi_{i,s} 
\]
We now estimate the $s$-derivative of $\log(F_{i,s})$. In the tubular neighborhood $U$, we can apply the formula (\ref{eq:reeb_vector_field}) for the Reeb vector field in Remark \ref{rmk:coordinate_formulas} to see that
\[
|\frac{d\alpha_{i,s}}{ds}(R_{i,s})| \le 2|\frac{d\alpha_{\chi_{i,s}}}{ds}(R_{\chi_{i,s}})| \le 2|\chi_1 - \chi_2| \cdot |d\theta(R_{\chi_{i,s}})| \le 2|d\theta(R_{\chi_{i,s}})| = \frac{2|\rho^*d\alpha(\partial_t,\partial_r)|}{G + H_{\chi_{i,s}}}
\]
Moreover, by (\ref{eq:cofinal_pa_1}) and the formula for the volume in (\ref{eq:volume_form}), we have $2(G + H_{\chi_{i,s}}) \ge G$ and thus
\begin{equation} \label{eq:alpha_d_estimate}
 |\frac{d\alpha_{i,s}}{ds}(R_{i,s})| \le \frac{2|\rho^*d\alpha(\partial_t,\partial_r)|}{G + H_{\chi_{i,s}}} \le \frac{4|\rho^*d\alpha(\partial_t,\partial_r)|}{G} = 4|\rho^*d\theta(R)| \le \log(C)
\end{equation}
On the other hand, the estimate (\ref{eq:alpha_d_estimate}) holds outside of $U$ since $\alpha_{i,s}$ is $s$-independent (and equal to $\rho^*\alpha$) in that part of $Y$. By integrating this inequality, we get the estimate
\[C^{-1} \le |F_{i,1}| \le C\]
If we let $\phi_i = \phi_{i,1}$ and $\beta = \alpha_1$, then this establishes strong cofinality.
\[
 C^{-1} \beta \le \phi_i^*\alpha_i \le C \beta \qquad L_{i+1} > D^2L_i \qquad D > C > 1 \qedhere
\]\end{proof}

\noindent As a corollary of Lemmas \ref{lem:universally_tight}, \ref{lem:tracking_hypertight} and \ref{lem:cofinal_for_pseudoanosov}, we get Theorem \ref{thm:hypertight} from the introduction.

\begin{corollary} A pseudo-Anosov contact structure is asymptotically hypertight and universally tight.
\end{corollary}

\subsection{Proofs Of Main Results.} In this part, we prove several of the results claimed in the introduction. We start by giving the relevant computation of the cylindrical contact homology of pseudo-Anosov contact manifolds.

\begin{theorem}[Theorem \ref{thm:contact_homology}] \label{thm:body_contact_homology} Let $\xi$ be a pseudo-Anosov contact structure on a closed 3-manifold $Y$. Then for any primitive free homotopy class $\Gamma$ and any pseudo-Anosov Reeb flow $\Phi$ inducing $\xi$, we have
\[
CH(Y,\xi;\Gamma) \neq 0 \qquad\text{if and only if}\qquad \Gamma \in \mathcal{P}(\Phi)
\]
Moreover, for any contact form $\alpha$ on $(Y,\xi)$, there is a constant $C > 0$ such that if $\Phi$ has a closed orbit in the homotopy class $\Gamma$ of period $L$ and the homotopy class $\Gamma$ contains no singular orbits, then
\[
CH^{CL}(Y,\alpha;\Gamma) \neq 0
\]
\end{theorem}

\begin{proof} We break the proof into several cases depending on the types of orbits present in the homotopy class $\Gamma$ (including the case where $\Gamma$ contains no orbits).

\vspace{3pt}

{\bf Case 1: Class Contains Orbits.} We start with the case where $\Gamma$ is present in $\mathcal{P}(\Phi)$. Let $A$ be any constant larger than the maximal period of all simple singular orbits of $\Phi$. By Lemma \ref{lem:cofinal_for_pseudoanosov}, there exists a sequence of smooth non-degenerate contact forms $\alpha_i$, actions $L_i$ that are not in the spectrum of $\Phi$ and contactomorphisms isotopic to the identity
\[\phi_i:(Y,\xi) \to (Y,\on{ker}(\alpha_{\chi_i}))\] such that the contact forms $\alpha_{\chi_i}$ satisfy
\[
\text{$\alpha_i$ tracks the Reeb flow of $\alpha$ to actions $A$ and $L_i$} \qquad\text{and}\qquad \text{$(\phi^*_i\alpha_i,L_i)$ \text{ is strongly cofinal}}
\]
By Lemma \ref{lem:non_vanishing_CCH} and Remark \ref{rmk:grading_nonvanishing}, it suffices to show that there is a choice of $A$ such that, for all sufficiently large $i$ and for all $S$ with $A \le S < L_i$, the map
\begin{equation} \label{eq:nonvanishing_proof}
CH_S(Y,\alpha_i;\Gamma) \to CH_{L_i}(Y,\alpha_i;\Gamma) \qquad\text{is injective and non-zero in a fixed grading}
\end{equation}
We next break this case into three sub-cases based on the types of orbits present in the homotopy class $\Gamma$ according to Lemma \ref{lem:orbit_types_in_homotopy_classes}.

\vspace{3pt}

{\bf Case 1a: Positive Hyperbolic.} First assume that $\Gamma$ is represented by (possibly many) positive hyperbolic orbits and non-rotating singular orbits of $\Phi$. This is the most complicated case.

\vspace{3pt}

In this case, choose the constant $A$ to be larger than the maximal period of all the non-rotating singular orbits $\eta_1 \dots \eta_k$ in the homotopy class $\Gamma$, and greater than the period of at least one smooth orbit in the homotopy class $\Gamma$. We consider the chain complex for the contact homology of $(\alpha_i,L_i)$ in the class $\Gamma$ up to action $A \le S < L_i$, and the corresponding contact homology groups
\[
C_S(Y,\alpha_i;\Gamma) \qquad\text{and}\qquad CH_S(Y,\alpha_i;\Gamma) = H(C_S(Y,\alpha_i;\Gamma))
\]

First, we claim that $CH_S(Y,\alpha_i;\Gamma) \neq 0$. Note that since $\alpha_i$ is $L_i$-hypertight, this complex can be written as a direct sum
\[
C_S(Y,\alpha_i;\Gamma) = C^S_0 \oplus C_1^S
\]
where $C_i^S$ is generated by the orientation lines of closed orbits of grading $i \mod 2$ up to period $S$. By Example \ref{ex:Reeb_Lef}, the Euler characteristic of this complex is precisely given by
\[
\on{dim}(C^S_0) - \on{dim}(C^S_1) = \sum_{\gamma} \on{Lef}(\gamma) 
\]
where the sum is over all closed Reeb orbits $\gamma$ of $\gamma_i$ in the homotopy class $\Gamma$ of period less than $S$ and $\on{Lef}(\gamma)$ denotes the Lefschetz sign of the orbit. Since $\alpha_i$ is a tracking approximation of $\alpha$ up to action $S$, we can compute this sum of Lefschetz numbers as the sum of the corresponding Lefschetz numbers of smooth and singular orbits of $\Phi$ up to action $S$ in the class $\Gamma$. Since all such orbits are positive hyperbolic or non-rotating singular, Lemma \ref{lem:tracking_Lef} and Example \ref{ex:PA_Lef} then imply
\[
\sum_{\gamma} \on{Lef}(\gamma) = \sum_{\eta} p_\eta - 1 
\]
where the sum on the right is over closed orbits of the pseudo-Anosov Reeb flow $\Phi$ of period $L_i$ or less and $p_\eta$ denotes the number of prongs of $\eta$ (which is $2$ when $\eta$ is a smooth orbit). In particular, we have
\[
\chi(C_S(Y,\alpha_i;\Gamma)) = \on{dim}(C^S_0) - \on{dim}(C^S_1) \ge 1
\]
In particular, this implies that the homology $CH_S(Y,\alpha_i;\Gamma)$ is non-zero. 

\vspace{3pt}

Second, we claim that the map $CH_S(Y,\alpha_i;\Gamma) \to CH_{L_i}(Y,\alpha_i;\Gamma)$ is injective in grading zero. For this, we note that every orbit of action between $A \le S$ and $L_i$ is a (non-singular) positive hyperbolic orbit since every singular orbit has action less than $A$. In particular, these orbits have grading zero and thus
\[
C^A_1 = C^S_1 = C^{L_i}_1 \qquad\text{for all $S$ with $A \le A \le L_i$}
\]
This implies that the map $CH_S(Y,\alpha_i) \to CH_{L_i}(Y,\alpha_i)$ is injective on the zero graded part.

\vspace{3pt}

{\bf Case 1b: Negative Hyperbolic.} Next assume $\Gamma$ is represented by a unique (non-singular) negative hyperbolic orbit $\eta$. In this case, choose $A$ be any constant larger than the minimal period of $\eta$. Since $\alpha_i$ converges smoothly to $\alpha$ along $\eta$, the contact forms $\alpha_i$ have a unique hyperbolic orbit of period less than $A$ for all sufficiently large $i$. Since $\alpha_i$ tracks $\alpha$ up to action $L_i$, this is the unique orbit in the class $\Gamma$ up to action $L_i$. The injectivity of (\ref{eq:nonvanishing_proof}) follows.

\vspace{3pt}

{\bf Case 1c: Rotating Singular.} Next assume $\Gamma$ is represented by a unique rotating singular orbit $\eta$. In their case, choose $A$ to be any constant larger than the minimal period of $\eta$. Lemma \ref{lem:tracking_Lef} and Example \ref{ex:PA_Lef} state that the relative Lefschetz index of the local return map of this singular orbit is 1. Thus an analogous argument to Cases 1 and 2 shows (\ref{eq:nonvanishing_proof}). 

\vspace{3pt}

{\bf Case 2: No Orbits.} Finally, assume that $\Gamma$ contains no closed orbits. We again may choose a cofinal sequence $(\phi^*\alpha_i,L_i)$ with $\alpha_i$ tracking $\alpha$ up to action $L_i$. In this case, since every closed orbit of $\alpha_i$ of action bounded by $L_i$ is homotopic to an orbit of $\Phi$, we have
\[
CH_{L_i}(Y,\alpha_i;\Gamma) = 0 \qquad\text{and thus}\qquad CH(Y,\xi;\Gamma) = \underset{i}{\on{colim}}\big(CH_{L_i}(Y,\alpha_i;\Gamma)\big) = 0\]

{\bf Final Part.} For the final claim of the theorem, fix a contact form $\alpha$ on $(Y,\xi)$ and a free homotopy class $\Gamma$ containing a smooth closed orbit of period $L$. By scaling $\alpha$, we may assume that $\alpha \le \alpha_i$. Then Lemma \ref{lem:non_vanishing_CCH} and Remark \ref{rmk:grading_nonvanishing} state that
\[
CH^A(Y,\alpha;\Gamma) \neq 0
\]
where $A$ is the chosen constant in the cases above. Now simply note that, in Cases 1 and 2 when $\Gamma$ contains no singular orbits, we may choose $A$ to be any constant larger than $L$. In particular
\[
CH^{2L}(Y,\alpha;\Gamma) \neq 0 \qedhere
\]\end{proof}

\begin{theorem}[Theorem \ref{thm:entropy}] Let $\xi$ be a pseudo-Anosov contact structure. Then any smooth Reeb flow $\Phi$ for $\xi$ has positive entropy.
\end{theorem}

\begin{proof} Fix a pseudo-Anosov Reeb flow $\Psi$ on $(Y,\xi)$ and a contact form $\alpha$ generating a smooth Reeb flow $\Phi$. By Lemma \ref{lem:CHG_entropy_bound} and Lemma \ref{lem:CHG_properties}, it suffices to show that the CH homotopical growth rate $\on{CHGr}(Y,\alpha)$ is positive for any smooth contact form $\alpha$. By Theorem \ref{thm:body_contact_homology}, there is a constant $C > 0$ such that if $\Phi$ has a closed orbit in a free homotopy class $\Gamma$ of period $L$, then $\Psi$ has a closed orbit in the class $\Gamma$ of period $CL$ or less. It follows that
\[
\on{CHGr}(Y,\alpha) \ge C^{-1} \cdot \on{Gr}(Y,\Psi)
\]
The latter quantity is positive and equal to the entropy of $\Psi$ by Theorem \ref{thm:PA_entropy}.\end{proof}

\bibliographystyle{hplain}
\bibliography{standard_bib}

\end{document}